\documentclass[12pt,reqno,a4paper]{amsart}
\usepackage{amssymb,amscd,amsmath}
\usepackage{pspicture}
\usepackage{graphicx}

\usepackage[matrix,arrow,curve,frame]{xy}    

\xymatrixcolsep{1.9pc}                          
\xymatrixrowsep{1.9pc}
\newdir{ >}{{}*!/-5pt/\dir{>}}                  

 \calclayout

\makeatletter

\renewcommand{\subsection}{\@startsection{subsection}{1}{0pt}{-3.25ex plus -1ex minus-.2ex}{1.5ex plus.2ex}{\normalfont\it}}
\renewcommand{\section}{\@startsection{section}{1}{\parindent}{3.5ex plus 1ex minus .2ex}{2.3ex plus.2ex}{\sc}}

\renewcommand{\phi}{\varphi}
\renewcommand{\leq}{\leqslant}
\renewcommand{\geq}{\geqslant}
\renewcommand{\epsilon}{\varepsilon}

\renewcommand{\kappa}{\varkappa}

\DeclareMathOperator{\Ext}{Ext}
\DeclareMathOperator{\spec}{Spec}
 \DeclareMathOperator{\cyl}{cyl}
\DeclareMathOperator{\Nvmod}{Nvmod}

 \DeclareMathOperator{\Map}{Map}
 \DeclareMathOperator{\mot}{mot}

\DeclareMathOperator{\Hom}{Hom} 
 
 \DeclareMathOperator{\id}{id}

 \DeclareMathOperator{\Mor}{Mor}
 
\DeclareMathOperator{\Ho}{Ho}

\DeclareMathOperator{\coker}{Coker} \DeclareMathOperator{\nis}{nis}
 \DeclareMathOperator{\Ar}{Ar}
\DeclareMathOperator{\lm}{lim} 
 \DeclareMathOperator{\Mod}{Mod}
 \DeclareMathOperator{\Ob}{Ob}
\DeclareMathOperator{\supp}{Supp}

\newcommand{\lo}{\lm}

\newcommand{\lra}[1]{\bl{#1}\longrightarrow\relax}
\newcommand{\bl}[1]{\buildrel #1\over}
\newcommand{\cc}{\mathcal}
\newcommand{\bb}{\mathbb}
\newcommand{\ps}{\oplus}

\newcommand{\op}{{\textrm{\rm op}}}

\newcommand{\wh}{\widehat}
\newcommand{\wt}{\widetilde}

\newcommand{\cork}{Cor_K}
\newcommand{\corvirt}{Cor_{virt}}
\newcommand{\corkplus}{Cor_{K^\oplus}}
\newcommand{\ifff}{if and only if }

\newcommand{\sheff}{SH^{\mot}}
\newcommand{\shnis}{SH^{\nis}}

\newtheorem{thm}{Theorem}[section]
\newtheorem{prop}[thm]{Proposition}
\newtheorem*{sublem}{Sublemma}
\newtheorem{cor}[thm]{Corollary}
\newtheorem{lem}[thm]{Lemma}

\newtheorem{rem}[thm]{Remark}

\newtheorem{exs}[thm]{Examples}
\newtheorem{example}[thm]{Example}
\newtheorem{defs}[thm]{Definition}

\makeatother

\begin{document}

\footskip30pt


\title{K-motives of algebraic varieties}
\author{Grigory Garkusha}
\address{Department of Mathematics, Swansea University, Singleton Park, Swansea SA2 8PP, United Kingdom}
\email{g.garkusha@swansea.ac.uk}

\author{Ivan Panin}
\address{St. Petersburg Branch of V. A. Steklov Mathematical Institute,
Fontanka 27, 191023 St. Petersburg, Russia}
\email{paniniv@gmail.com}

\thanks{This paper was written during the visit of the second author to
Swansea University supported by EPSRC grant EP/H021566/1. He would
like to thank the University for the kind hospitality}

\keywords{Motivic homotopy theory, algebraic $K$-theory, spectral
categories}

\subjclass[2000]{14F42, 19E08, 55U35}

\begin{abstract}
A kind of motivic algebra of spectral categories and modules over
them is developed to introduce $K$-motives of algebraic varieties.
As an application, bivariant algebraic $K$-theory $K(X,Y)$ as well
as bivariant motivic cohomology groups $H^{p,q}(X,Y,\bb Z)$ are
defined and studied. We use Grayson's machinery~\cite{Gr} to produce
the Grayson motivic spectral sequence connecting bivariant
$K$-theory to bivariant motivic cohomology. It is shown that the
spectral sequence is naturally realized in the triangulated category
of $K$-motives constructed in the paper. It is also shown that
ordinary algebraic $K$-theory is represented by the $K$-motive of
the point.
\end{abstract}
\maketitle

\thispagestyle{empty} \pagestyle{plain}

\newdir{ >}{{}*!/-6pt/@{>}} 

\tableofcontents

\section{Introduction}

The triangulated category of motives $DM^{eff}$ in the sense of
Voevodsky~\cite{Voe1} does not provide a sufficient framework to
study such a fundamental object as the motivic spectral sequence. In
this paper we construct the triangulated category of $K$-motives
over any field $F$. This construction provides a natural framework
to study (bivariant) $K$-theory in the same fashion as the
triangulated category of motives provides a framework for motivic
cohomology.

The main idea is to use formalism of spectral categories over smooth
schemes $Sm/F$ and modules over them. In this language a transfer
from one scheme $X$ to another scheme $Y$ is a symmetric spectrum
$\cc O(X,Y)$ such that there is an associative composition law
   $$\cc O(Y,Z)\wedge\cc O(X,Y)\to\cc O(X,Z).$$
The category of $\cc O$-modules $\Mod\cc O$ consists of presheaves
of symmetric spectra having ``$\cc O$-transfers". The main spectral
categories we work with are $\cc O_{K^{Gr}}$, $\cc O_{K^\oplus}$,
$\cc O_K$. They come from various symmetric $K$-theory spectra
associated with the category of bimodules $\cc P(X,Y)$, $X,Y\in
Sm/F$. By a bimodule we mean a coherent $\cc O_{X\times Y}$-module
$P$ such that $\supp P$ is finite over $X$ and the coherent $\cc
O_X$-module $(p_X)_*(P)$ is locally free.

In order to develop a satisfactory homotopy theory of presheaves of
symmetric spectra with $\cc O$-transfers we specify the condition of
being a ``motivically excisive spectral category". In this case one
produces a compactly generated triangulated category $\sheff\cc O$
which plays the role of the triangulated category of motives
associated with the spectral category of transfers $\cc O$.
Voevodsky's category $DM^{eff}$ can be recovered in this way from
the Eilenberg--Mac~Lane spectral category $\cc O_{cor}$ associated
with the category of correspondences. The spectral categories $\cc
O_{K^{Gr}}$, $\cc O_{K^\oplus}$, $\cc O_K$ produce equivalences of
triangulated categories 
   $$\sheff\cc O_{K^{Gr}}\simeq\sheff\cc O_{K^\oplus}\simeq\sheff\cc O_K.$$

The $K$-motive $M_K(X)$ of a smooth scheme $X$ over $F$ is the image
of the free $\cc O_K$-module $\cc O_K(-,X)$ in $\sheff\cc O_K$ (see
Definition~\ref{refer}). Then there is an isomorphism (see
Corollary~\ref{aqemm})
   $$K_i(X)\cong\sheff\cc O_K(M_K(X)[i],M_K(pt)),\quad i\in\bb Z.$$
Thus ordinary $K$-theory is represented by the $K$-motive of the
point.

One of the main computational tools of the paper is the ``Grayson
motivic spectral sequence". It is a strongly convergent spectral
sequence of the form
   $$E_2^{pq}=H_{\cc M}^{p-q,-q}(X,\bb Z)\Longrightarrow K_{-p-q}(X)$$
where the groups on the left hand side are motivic cohomology groups
of $X$. We show in Theorem~\ref{spektralka3} that it is recovered
from the ``Grayson tower" in $\sheff\cc O_{K^{Gr}}$
   $$\cdots\xrightarrow{f_{q+1}}M_{K^{Gr}}(q)(pt)\xrightarrow{f_q}M_{K^{Gr}}(q-1)(pt)\xrightarrow{f_{q-1}}\cdots\xrightarrow{f_1}M_{K^{Gr}}(pt),$$
where $M_{K^{Gr}}(q)(pt)$-s are certain $\cc O_{K^{Gr}}$-modules.
Thus the triangulated category of $K$-motives provides a sufficient
framework to study the motivic spectral sequence. In fact, we
construct its bivariant counterpart. It is used to show that $\cc
O_{K^{Gr}}$, $\cc O_{K^\oplus}$, $\cc O_K$ are motivically excisive
spectral categories.

One should stress that we do not construct a tensor product on
$\sheff\cc O_K$. 

Throughout the paper we denote by $Sm/F$ the category of smooth
separated schemes of finite type over the base field $F$.

\section{Model structures for symmetric spectra}\label{putrya}

In this section we collect basic facts about symmetric spectra of
simplicial sets $Sp^\Sigma$ we shall need later. We refer the reader
to~\cite{HSS,Sch} for details.

\begin{defs}{\rm
(1) An object $A$ of a model category $\cc M$ is {\it finitely
presentable\/} if the set-valued Hom-functor $\Hom_{\cc M}(A,-)$
commutes with all filtered colimits.

(2) Following~\cite{DRO} a cofibrantly generated model category $\cc
M$ is {\it weakly finitely generated\/} if $I$ and $J$ can be chosen
such that the following conditions hold.

\begin{itemize}
\item[$\diamond$] The domains and the codomains of the maps in $I$ are finitely
presentable.

\item[$\diamond$] The domains of the maps in $J$ are small.

\item[$\diamond$] There exists a subset $J'$ of $J$ of maps with finitely presentable
domains and codomains, such that a map $f:A\to B$ in $\cc M$ with
fibrant codomain $B$ is a fibration if and only if it is contained
in $J'-inj$.
\end{itemize}

}\end{defs}

Recall from \cite[Chapter V]{H} that it is possible to define
$X\otimes K$ for an object $X$ in a model category and a simplicial
set $K$, even if the model category is not simplicial.

\begin{lem}\label{pqp}
Let $\cc M$ be a left proper, cellular, weakly finitely generated
model category in which all objects are small and let $S$ be a set
of cofibrations in $\cc M$. Suppose that, for every domain or
codomain $X$ of $S$ and every finite simplicial set $K$, $X\otimes
K$ is finitely presentable. Then the Bousfield localization $\cc
M/S$ is weakly finitely generated.
\end{lem}

\begin{proof}
The proof is like that of~\cite[4.2]{Hov}.
\end{proof}

We first define level projective and flat model structures for
symmetric spectra. A morphism $f:A\to B$ of symmetric spectra is
called a {\it flat cofibration\/} if and only if for every level
cofibration $g:X\to Y$ the pushout product map
   $$f\wedge g:B\wedge X\cup_{A\wedge X}A\wedge Y\to B\wedge Y$$
is a level cofibration. In particular, every flat cofibration is a
level cofibration.

\begin{thm}[\cite{HSS,Sch}]\label{spect}
The category of symmetric spectra $Sp^\Sigma$ admits the following
two level model structures in which the weak equivalences are those
morphisms $f:X\to Y$ such that for all $n\geq 0$ the map $f_n:X_n\to
Y_n$ is a weak equivalence of simplicial sets.

$(1)$ In the projective level model structure a morphism $f:X\to Y$
is a projective level fibration if and only if for every $n\geq 0$
the map $f_n:Y_n\to X_n$ is a Kan fibration of simplicial sets. A
morphism $f:X\to Y$ is a {\it projective cofibration\/} if it has
the left lifting property with respect to all acyclic projective
level fibrations.

$(2)$ In the flat level model structure the cofibrations are the
flat cofibrations. A morphism $f:X\to Y$ is a {\it flat level
fibration\/} if it has the right lifting property with respect to
all acyclic flat cofibrations.

The two stable model structures are cellular, proper, simplicial,
weakly finitely generated, symmetric monoidal with respect to the
smash product of symmetric spectra and satisfy the monoid axiom in
the sense of Schwede--Shipley~\cite{SS1}.
\end{thm}

\begin{proof}
By~\cite{HSS,Sch} (1) and (2) determine proper, simplicial,
cofibrantly generated and symmetric monoidal model structures.
By~\cite[3.2.13]{HSS} every symmetric spectrum is small.
By~\cite{HSS,Sch} the domains and codomains of generating (trivial)
cofibrations are finitely presentable and projective (flat)
cofibrations are levelwise injections of simplicial sets. It follows
that both model structures are cellular and weakly finitely
generated. The monoid axiom follows from~\cite[III.1.11]{Sch} and
the fact that in a weakly finitely generated model category the
class of weak equivalences is closed under filtered colimits
by~\cite[3.5]{DRO}.
\end{proof}

Recall that a symmetric spectrum $X$ is {\it injective\/} if for
every monomorphism which is also a level equivalence $i:A\to B$ and
every morphism $f:A\to X$ there exists an extension $g:B\to X$ with
$f=gi$. A morphism $f:A\to B$ of symmetric spectra is a {\it stable
equivalence\/} if for every injective $\Omega$-spectrum $X$ the
induced map $[f,X]:[B,X]\to [A,X]$ on homotopy classes of spectrum
morphisms is a bijection.

\begin{thm}[\cite{Hov,HSS,Sch}]\label{spectra}
The category of symmetric spectra $Sp^\Sigma$ admits the following
two stable model structures in which the weak equivalences are the
stable equivalences.

$(1)$ In the projective stable model structure the cofibrations are
the projective cofibrations. A morphism $f:X\to Y$ is a {\it stable
projective fibration\/} if it has the right lifting property with
respect to all acyclic projective cofibrations.

$(2)$ In the flat stable model structure the cofibrations are the
flat cofibrations. A morphism $f:X\to Y$ is a {\it stable flat
fibration\/} if it has the right lifting property with respect to
all acyclic flat cofibrations.

The two stable model structures are cellular, proper, simplicial,
weakly finitely generated, symmetric monoidal with respect to the
smash product of symmetric spectra and satisfy the monoid axiom in
the sense of Schwede--Shipley~\cite{SS1}.
\end{thm}

\begin{proof}
By~\cite{HSS,Sch} (1) and (2) determine proper, simplicial,
cofibrantly generated and symmetric monoidal model structures. The
domains and codomains of generating cofibrations are finitely
presentable and by~\cite[3.2.13]{HSS} every symmetric spectrum is
small. By~\cite{HSS,Sch} projective (flat) cofibrations are
levelwise injections of simplicial sets. It follows that both model
structures are cellular. By~\cite[III.2.2]{Sch} and~\cite[1.1.11]{H}
stable equivalences which are monomorphisms of symmetric spectra are
closed under pushouts. By~\cite[5.4.1]{HSS} the monoid axiom is true
for the stable projective model structure. The monoid axiom for the
stable flat model structure is proved like that of~\cite[5.4.1]{HSS}
(for this use as well~\cite[III.1.11]{Sch} and the fact that in a
weakly finitely generated model category the class of weak
equivalences is closed under filtered colimits by~\cite[3.5]{DRO}).

It remains to verify that both model structures are weakly finitely
generated. By~\cite[p.~109]{Hov} the stable projective model
structure on $Sp^\Sigma$ is the Bousfield localization of the level
projective model structure with respect to the set $\cc S$
   $$F_{n+1}(C\wedge S^1)\to F_n(C)$$
as $C$ runs through the domains and codomains of the generating
cofibrations of pointed simplicial sets $SSets_*$ and each $F_n$ is
the left adjoint to the $n$th evaluation functor $Ev_n:Sp^\Sigma\to
SSets_*$. It follows from Lemma~\ref{pqp} and Theorem~\ref{spect}
that the stable projective model structure on $Sp^\Sigma$ is weakly
finitely generated.

If we show that the stable flat model structure is the Bousfield
localization of the level flat model structure with respect to the
set $\cc S$, it will follow from Lemma~\ref{pqp} and
Theorem~\ref{spect} that it is weakly finitely generated. The
cofibrations in both model structures are the same. Since level
projective model structure is Quillen equivalent to level flat model
structure, then the corresponding Bousfield localizations with
respect to $\cc S$ are Quillen equivalent as well. Therefore weak
equivalences in the stable flat model structure and in the Bousfield
localization of the level flat model structure with respect to the
set $\cc S$ coincide as was to be shown.
\end{proof}

We want to make several remarks about symmetric spectra. One of
tricky points when working with symmetric spectra is that the stable
equivalences of symmetric spectra can not be defined by means of
stable homotopy groups. It is not enough to invert
$\pi_*$-isomorphisms (=stable weak equivalence of ordinary spectra)
to get a satisfactory homotopy category of symmetric spectra.
Instead one inverts a bigger class -- that of stable weak
equivalences between symmetric spectra in the sense of~\cite{HSS}.

Given a symmetric spectrum $X$ and its stably fibrant model $\gamma
X$ in $Sp^\Sigma$, stable homotopy groups of $X$ can considerably be
different from those of $\gamma X$ in general (see,
e.g.~\cite{HSS,Sch}). In particular, Hom-sets $\Ho(Sp)(S^i,X)$ are
different from $\Ho(Sp^\Sigma)(S^i,X)$. Nevertheless there is an
important class of {\it semistable\/} symmetric spectra within which
stable equivalences coincide with $\pi_*$-isomorphisms. Recall that
a symmetric spectrum is semistable if some (hence any) stably
fibrant replacement is a $\pi_*$-isomorphism. Here a stably fibrant
replacement is a stable equivalence $X\to\gamma X$ with target an
$\Omega$-spectrum.

Suspension spectra, Eilenberg--Mac Lane spectra, $\Omega$-spectra or
$\Omega$-spectra from some point $X_n$ on are examples of semistable
symmetric spectra (see~\cite[Example~I.4.48]{Sch}). So Waldhausen's
algebraic $K$-theory symmetric spectrum we shall discuss later is
semistable. Semistability is preserved under suspension, loop,
wedges and shift~\cite[Example~I.4.51]{Sch}.

In what follows we shall use these facts without further comments.

\section{Spectral categories and modules over them}

To define $K$-motives, we work in the framework of spectral
categories and modules over them in the sense of
Schwede--Shipley~\cite{SS}. We start with preparations.

A biexact functor of Walddhausen categories is a functor $\cc
A\times\cc B\to\cc C$, $(A,B)\mapsto A\otimes B$, having the
property that for every $A\in\cc A$ and $B\in\cc B$ the partial
functors $A\otimes-$ and $-\otimes B$ are exact, and for every pair
of cofibrations $A\rightarrowtail A'$ and $B\rightarrowtail B'$ in
$\cc A$ and $\cc B$ respectively, the map $A'\otimes
B\coprod_{A\otimes B}A\otimes B'\to B\otimes B'$ is a cofibration in
$\cc C$.

For example let $\cc P(X,Y)$, $X,Y\in Sm/F$, be the category of
coherent $\cc O_{X\times Y}$-modules $P$ such that $\supp P$ is
finite over $X$ and the coherent $\cc O_X$-module $(p_X)_*(P)$ is
locally free. For example, let $f:X\to Y$ be a morphism of smooth
schemes, and let $\Gamma_f$ be its graph; then $\Gamma_f\in\cc
P(X,Y)$. Let $X,Y,U\in Sm/F$ be three smooth schemes. In this case
we have a natural functor
   \begin{equation}\label{pxu}
    \cc P(X,Y)\times\cc P(Y,U)\to\cc P(X,U),\quad P\times Q\mapsto P\otimes Q:=(p_{X,U})_*(p^*_{X,Y}(P)\otimes_{\cc O_{X\times Y\times U}}p^*_{Y,U}(Q))
   \end{equation}
By~\cite[section~1]{Sus} the sheaf on the right really belongs to
$\cc P(X,U)$ and the above functor is biexact.

Let $\corvirt$ be a category whose objects are those of $Sm/F$ and
satisfying the following conditions:

\begin{itemize}
\item[(a)] For any pair of smooth schemes $X,Y\in Sm/F$ there is a Waldhausen category
$(\cc C_{virt}(X,Y),w)$ with $w$ a family of weak equivalences in
$\cc C_{virt}(X,Y)$ such that $\Ob\cc
C_{virt}(X,Y)=\Mor{\corvirt}(X,Y)$.

\item[(b)] For any triple $X,Y,Z\in Sm/F$ there is an associative biexact functor of
Waldhausen categories
   \begin{equation}\label{prod}
    \cc C_{virt}(X,Y)\times\cc C_{virt}(Y,Z)\lra{\phi}\cc C_{virt}(X,Z).
   \end{equation}
By associativity we mean that the diagram
   $$\xymatrix{\cc C_{virt}(X,Y)\times(\cc C_{virt}(Y,Z)\times\cc
   C_{virt}(Z,W))\ar[r]^(.58){1\times\phi}\ar[d]_\cong
               &\cc C_{virt}(X,Y)\times\cc C_{virt}(Y,W)\ar[dd]^{\phi}\\
               (\cc C_{virt}(X,Y)\times\cc C_{virt}(Y,Z))\times\cc C_{virt}(Z,W)\ar[d]_{\phi\times
               1}\\
               \cc C_{virt}(X,Z)\times\cc C_{virt}(Z,W)\ar[r]^{\phi}&\cc C_{virt}(X,W)}$$
is commutative.

\item[(c)] For any $X\in Sm/F$ there is a distinguished object $1_X\in\cc
C_{virt}(X,X)$ (the ``tensor product unit object") such that
   $$\phi(1_X,P)=\phi(P,1_Y)$$
for all $P\in\cc C_{virt}(X,Y)$.
\end{itemize}
In other words, the category $\corvirt$ is ``enriched" over
Waldhausen categories. Sometimes we refer to $\corvirt$ as the
category of ``virtual" correspondences.

\begin{example}{\rm
Given $X,Y\in Sm/F$, consider the category $\cc P(X,Y)$. It is an
exact category, and therefore it can be regarded as a Waldhausen
category with the family of weak equivalences being isomorphisms.
Clearly, biproduct~\eqref{pxu} is associative up to isomorphism but
not strictly associative. We replace the exact categories $\cc
P(X,Y)$-s by equivalent exact categories $\cc P'(X,Y)$-s and define
a strictly associative biproduct
   $$\phi=\phi_{XYU}:\cc P'(X,Y)\times\cc P'(Y,U)\to\cc P'(X,U).$$

One approach is to define $\cc P'(X,Y)$ as follows.
Consider sequences $(P_1,\ldots,P_k)$, $k\geq 1$, of objects in
other categories $\cc P(U_i,V_i)$, so that $U_1=X$, $V_k=Y$ and
$V_i=U_{i+1}$ for any $i<k$. If $X=Y$ the object $1_X=\Gamma_1$ of
$\cc P(X,X)$ is a monoidal unit. There is a relation for sequences
$(P)=(1_X,P)=(P,1_Y)$ with $P\in\cc P(X,Y)$. We also require
$(\Gamma_f,\Gamma_g)=(\Gamma_{gf})$ to hold for any morphisms of
smooth schemes $f:X\to U,g:U\to Y$.

We say that a sequence $(P_1,...,P_k)$ is of minimal length if it
can not be reduced to a sequence of smaller length by means of the
relations above. We set
   $$P_1\otimes\cdots\otimes P_k:=((\cdots(P_1\otimes P_2)\otimes\cdots)\otimes P_{k-1})\otimes P_k.$$
By definition, objects of $\cc P'(X,Y)$ are the sequences
$(P_1,\ldots,P_k)$ of minimal length. Define the arrows between two
sequences of minimal length $(P_1,\ldots,P_k)$ and
$(Q_1,\ldots,Q_l)$ by
   $$\Hom_{\cc P'(X,Y)}((P_1,...,P_k),(Q_1,...,Q_l)):=\Hom_{\cc P(X,Y)}(P_1\otimes\cdots\otimes P_k,Q_1\otimes\cdots\otimes Q_l).$$

One easily sees that $\cc P'(X,Y)$ is an exact category and the
natural exact functor $\cc P(X,Y)\to\cc P'(X,Y)$ sending an object
$P$ to the sequence of length one is an equivalence.

The tensor product
   $$\phi=\phi_{XYU}:\cc P'(X,Y)\times\cc P'(Y,U)\to\cc P'(X,U)$$
on the new objects is simply concatenation of sequences, which is
strictly associative. By construction,
   $$\phi(1_X,(P_1,...,P_k))=\phi((P_1,...,P_k),1_Y)=(P_1,...,P_k)$$
for any $(P_1,...,P_k)\in\cc P'(X,Y)$.

Let $\cork$ (respectively $\corkplus$) be the category whose objects
are those of $Sm/F$, $\Mor\cork(X,Y)=\Ob\cc P'(X,Y)$ (respectively
$\Mor\corkplus(X,Y)=\Ob\cc P'(X,Y)$) and for any $X,Y\in Sm/F$ let
$\cc C_{K}(X,Y)$ (respectively $\cc C_{K^\oplus}(X,Y)$) be the exact
category $\cc P'(X,Y)$ (respectively the same category $\cc P'(X,Y)$
but considered as an additive category). Note that the map taking a
morphism of smooth schemes $f:X\to Y$ to $\Gamma_f$ determines a
functor $Sm/F\to\cork$ (respectively a functor $Sm/F\to\corkplus$).

Using coherence properties for tensor product
(see~\cite[Ch.~VII]{Mac}), we have that for any triple $X,Y,Z\in
Sm/F$ there are biexact functors of Waldhausen categories
   $$\cc C_{K}(X,Y)\times\cc C_{K}(Y,Z)\lra{\phi}\cc C_{K}(X,Z)$$
and
   $$\cc C_{K^\oplus}(X,Y)\times\cc C_{K^\oplus}(Y,Z)\lra{\phi}\cc C_{K^\oplus}(X,Z)$$
satisfying conditions (b)-(c) above.

Thus $\cc P'$ yields two examples for $\cc C_{virt}$: $\cc C_{K}$
and $\cc C_{K^\oplus}$ respectively. We also get two categories of
correspondences $\cork$ and $\corkplus$ which are the same as
categories on smooth schemes but with different Waldhausen
categories on objects.

Given a (multisimplicial) additive category $\cc M$, we shall
sometimes write $K^\oplus\cc M$ to denote the $K$-theory symmetric
spectrum spectrum of $\cc M$.

}\end{example}

\begin{defs}\label{avtobus}{\rm
(1) Following~\cite{SS} a {\it spectral category\/} is a category
$\cc O$ which is enriched over the category $Sp^\Sigma$ of symmetric
spectra (with respect to smash product, i.e., the monoidal closed
structure of \cite[2.2.10]{HSS}). In other words, for every pair of
objects $o,o'\in\cc O$ there is a morphism symmetric spectrum $\cc
O(o,o')$, for every object $o$ of $O$ there is a map from the sphere
spectrum $S$ to $\cc O(o,o)$ (the ``identity element" of $o$), and
for each triple of objects there is an associative and unital
composition map of symmetric spectra $\cc O(o',o'')\wedge\cc O(o,o')
\to\cc O(o,o'')$. An $\cc O$-module $M$ is a contravariant spectral
functor to the category $Sp^\Sigma$ of symmetric spectra, i.e., a
symmetric spectrum $M(o)$ for each object of $\cc O$ together with
coherently associative and unital maps of symmetric spectra
$M(o)\wedge\cc O(o',o)\to M(o')$ for pairs of objects $o,o'\in\cc
O$. A morphism of $\cc O$-modules $M\to N$ consists of maps of
symmetric spectra $M(o)\to N(o)$ strictly compatible with the action
of $\cc O$. The category of $\cc O$-modules will be denoted by
$\Mod\cc O$.

(2) A {\it spectral functor\/} or a {\it spectral homomorphism\/}
$F$ from a spectral category $\cc O$ to a spectral category $\cc O'$
is an assignment from $\Ob\cc O$ to $\Ob\cc O'$ together with
morphisms $\cc O(a,b)\to\cc O'(F(a),F(b))$ in $Sp^\Sigma$ which
preserve composition and identities.

(3) The {\it monoidal product\/} $\cc O\wedge\cc O'$ of two spectral
categories $\cc O$ and $\cc O'$ is the spectral category where
$\Ob(\cc O\wedge\cc O'):=\Ob\cc O\times\Ob\cc O'$ and $\cc
O\wedge\cc O'((a,x),(b,y)):= \cc O(a,b)\wedge\cc O'(x,y)$.

(4) A {\it monoidal spectral category\/} consists of a spectral
category $\cc O$ equipped with a spectral functor $\diamond:\cc
O\wedge\cc O\to\cc O$, a unit $u\in\Ob\cc O$, a $Sp^\Sigma$-natural
associativity isomorphism and two $Sp^\Sigma$-natural unit
isomorphisms. Symmetric monoidal spectral categories are defined
similarly.

}\end{defs}

\begin{exs}\label{urr}{\rm
(1) A naive spectral category $\cc O_{naive}$ on $Sm/F$ is defined
as follows. $\cc O_{naive}$ has the same set of objects as $Sm/F$
and the morphism spectra are defined by
   $$\cc O_{naive}(X,Y)_p=\Hom_{Sm/F}(X,Y)_+\wedge S^p.$$
Here $S^p$ denotes the pointed simplicial set $S^p=
S^1\wedge\cdots\wedge S^1$ ($p$ factors) and the symmetric group
permutes the factors. Composition is given by the composite
   \begin{align*}
    \cc O_{naive}(Y,Z)_p\wedge\cc O_{naive}(X,Y)_q&=(\Hom_{Sm/F}(Y,Z)_+\wedge S^p)\wedge(\Hom_{Sm/F}(X,Y)_+\wedge S^q)\\
    &\xrightarrow{shuffle}\Hom_{Sm/F}(Y,Z)_+\wedge\Hom_{Sm/F}(X,Y)_+\wedge S^p\wedge S^q\\
    &\to\Hom_{Sm/F}(X,Z)_+\wedge S^{p+q}=\cc O_{naive}(X,Z)_{p+q}.
   \end{align*}
$\cc O_{naive}$ is a symmetric monoidal spectral category equipped
with a spectral functor $\cc O_{naive}\wedge\cc O_{naive}\to\cc
O_{naive}$, $\cc O_{naive}(X,Y)\wedge\cc O_{naive}(U,V)\to\cc
O_{naive}(X\times U,Y\times V)$, and a unit $\spec F\in\Ob\cc
O_{naive}$. It is straightforward to verify that the category of
$\cc O_{naive}$-modules can be regarded as the category of
presheaves $Pre^\Sigma(Sm/F)$ of symmetric spectra on $Sm/F$. This
is used in the sequel without further comment.

(2) Any {\it ringoid\/} $\cc A$, that is a category whose Hom-sets
are abelian groups with bilinear composition, gives rise to a
spectral category $H\cc A$ also called the Eilenberg--Mac Lane
spectral category of $\cc A$. In more detail, $H\cc A$ has the same
set of objects as $\cc A$ and the morphism spectra are defined by
$H\cc A(a,b)_p =\cc A(a,b)\otimes\wt{\bb Z}[S^p]$. Here $\wt{\bb
Z}[S^p]$ denotes the reduced simplicial free abelian group generated
by the pointed simplicial set $S^p$ and the symmetric group permutes
the factors. Composition is defined as above.

In Appendix we shall present another sort of Eilenberg--Mac Lane
spectral categories associated with ringoids. It will always be
clear from the context which of these sorts is used.

An important example of a ringoid is the category correspondences.
For any $X,Y\in Sm/F$ define $Cor(X,Y)$ to be the free abelian group
generated by closed integral subschemes $Z\subset X\times_FY$ which
are finite and surjective over a component of $X$. Let $X,Y,W\in
Sm/F$ be smooth schemes and let $Z\in Cor(X,Y),T\in Cor(Y,W)$ be
cycles on $X\times Y$ and $Y\times W$ each component of which is
finite and surjective over a component of $X$ (respectively over a
component of $Y$). One checks easily that the cycles $Z\times W$ and
$X\times T$ intersect properly on $X\times Y\times W$ and each
component of the intersection cycle $(Z\times W)_\bullet(X\times T)$
is finite and surjective over a component of $X$. Thus setting
$T\circ Z=(pr_{1,3})_*((Z\times Y)_\bullet(X\times T))$ we get a
bilinear composition map
   $$Cor(Y,W)\times Cor(X,Y)\to Cor(X,W).$$
In this way we get a ringoid (denoted $SmCor/F$) whose objects are
those of $Sm/F$ and $\Hom_{SmCor/F}(X,Y)=Cor(X,Y)$ ---
see~\cite{Voe} for details. The Eilenberg--Mac~Lane spectral
category corresponding to the ringoid will be denoted by $\cc
O_{cor}$.

For any smooth schemes $X,Y,U,V$ the external product of cycles
defines a homomorphism
   $$Cor(X,Y)\otimes Cor(U,V)\to Cor(X\times U,Y\times V)$$
which gives the structure of symmetric monoidal spectral category
for $\cc O_{cor}$\label{symmcor}.

(3) There are two other important ringoids $K_0^\oplus$ and $K_0$ on
$Sm/F$. Namely, for any $X,Y\in Sm/F$ define $K_0^\oplus(X,Y)$
(respectively $K_0(X,Y)$) to be the abelian group for the split
exact (respectively exact) category $\cc P'(X,Y)$. Composition is
given by tensor product. Denote by $\cc O_{K_0^\oplus}$ and $\cc
O_{K_0}$ their Eilenberg--Mac Lane spectral categories. There are
canonical homomorphisms
   $$K_0^\oplus(X,Y)\to K_0(X,Y)\to Cor(X,Y),\quad X,Y\in Sm/F.$$
Here, the first map is the obvious surjective homomorphism. The
second one takes the class $[P]$ of the coherent sheaf $P\in\cc
P(X,Y)$ to $\sum_Z\ell_{\cc O_{X\times Y,z}}P_z\cdot[Z]$, where the
sum is taken over all closed integral subschemes $Z\subset X\times
Y$ that are finite and surjective over a component of $X$ and $z$
denotes the generic point of the corresponding scheme $Z$.

By~\cite[section~1]{Sus} and~\cite[section~6]{Wlk} the canonical
homomorphisms yield maps between ringoids
   $$K_0^\oplus\to K_0\to Cor.$$
These induce spectral functors between spectral categories
   $$\cc O_{K_0^\oplus}\to\cc O_{K_0}\to\cc O_{cor}.$$
We shall prove below that the spectral categories $\cc
O_{K_0^\oplus}$ and $\cc O_{K_0}$ are symmetric monoidal (see
Corollary~\ref{qoqoqo}).

}\end{exs}

We want to construct a spectral category out of the category of
virtual correspondences $\corvirt$. Let $X,Y\in Sm/F$ and let $\cc
O_{virt}(X,Y):=K(\cc C_{virt}(X,Y))$ be the $K$-theory spectrum of
the Waldhausen category $\cc C_{virt}(X,Y)$. By~\cite[6.1.1]{GH}
and~\cite[Example~I.2.11]{Sch} $\cc O_{virt}(X,Y)$ naturally has the
structure of a symmetric spectrum. Note that $\cc
O_{virt}(X,Y)_0=|N.w\cc C_{virt}(X,Y)|$.

It follows from~\cite[section~6.1]{GH} and~\cite[Example~2.11]{Sch}
that associative biexact functors~\eqref{prod} induce an associative
law
   $$\cc O_{virt}(Y,Z)\wedge\cc O_{virt}(X,Y)\to\cc O_{virt}(X,Z).$$
Moreover, for any $X\in Sm/F$ there is a map $\mathbf 1:S\to\cc
O_{virt}(X,X)$ which is subject to the unit coherence law
(see~\cite[section~6.1]{GH}). Note that $\mathbf 1_0:S^0\to\cc
O_{virt}(X,X)_0$ is the map which sends the basepoint to the null
object and the non-basepoint to the unit object $1_X$ for the tensor
product. Thus the triple $(\cc O_{virt},\wedge,\mathbf 1)$
determines a spectral category on $Sm/F$.

Observe that $\cc O_{virt}$ is a symmetric monoidal spectral
category provided that there is a biexact functor
   $$\cc C_{virt}(X,Y)\times\cc C_{virt}(U,V)\to\cc C_{virt}(X\times U,Y\times V)$$
for all $X,Y,U,V\in Sm/F$ satisfying natural associativity,
symmetry, unit isomorphisms (where $\spec F$ is a unit).

In what follows the category of $\cc O_{virt}$-modules will be
denoted by $\cc M_{virt}$.

\begin{prop}\label{zdorovo}
The map
   $$\corvirt\to\cc M_{virt},\quad X\mapsto\cc O_{virt}(-,X),$$
determines a fully faithful functor.
\end{prop}

\begin{proof}
It follows from~\cite[sections~2.1-2.2]{DRO} that
   $$\Hom_{\cc M_{virt}}(M,N)=\Hom_{Sp^\Sigma}(S,{Sp}^\Sigma(M,N)).$$
Using ``Enriched Yoneda Lemma" one has,
   \begin{gather*}
    \Hom_{Sp^\Sigma}(S,{Sp}^\Sigma(\cc O_{virt}(-,X),\cc O_{virt}(-,Y)))\cong\Hom_{Sp^\Sigma}(S,\cc O_{virt}(X,Y))=\\
    \Hom_{SSets}(\Delta^0,|N.w\cc C_{virt}(X,Y)|)=|N.w\cc C_{virt}(X,Y)|_0=\corvirt(X,Y).
   \end{gather*}
Our statement now follows.
\end{proof}

If $\corvirt$ is either $\cork$ or $\corkplus$ then we shall denote
the corresponding spectral categories by $\cc O_K$ and $\cc
O_{K^\oplus}$ respectively. There is another spectral category $\cc
O_{K^{Gr}}$ we shall use later associated with $\corkplus$. It is
equivalent to $\cc O_{K^\oplus}$ and is based on
$S^\oplus$-construction of Grayson~\cite{Gr}. We let $Ord$ denote
the category of finite nonempty ordered sets. For $A\in Ord$ we
define a category $Sub(A)$ whose objects are the pairs $(i,j)$ with
$i\leq j\in A$, and where there is an (unique) arrow
$(i',j')\to(i,j)$ exactly when $i'\leq i\leq j\leq j'$. Given an
additive category $\cc M$, we say that a functor $M:Sub(A)\to\cc M$
is {\it additive\/} if $M(i,i)=0$ for all $i\in A$, and for all
$i\leq j\leq k\in A$ the map $M(i,k)\to M(i,j)\oplus M(j,k)$ is an
isomorphism. Here 0 denotes a previously chosen zero object of $\cc
M$. The set of such additive functors is denoted by $Add(Sub(A),\cc
M)$.

We define the simplicial set $S^\oplus\cc M$ by setting
   $$(S^\oplus\cc M)(A)=Add(Sub(A),\cc M).$$
An $n$-simplex $M\in S_n^\oplus\cc M$ may be thought of as a
compatible collection of direct sum diagrams $M(i,j)\cong
M(i,i+1)\oplus\cdots\oplus M(j-1,j)$. There is a natural map
$S^\oplus\cc M\to S\cc M$ which converts each direct sum diagram
$M(i,k)\cong M(i,j)\oplus M(j,k)$ into the short exact sequence
$0\to M(i,j)\to M(i,k)\to M(j,k)\to 0$.

Given a finite set $Q$, one can define the $|Q|$-fold iterated
$S^\oplus$-construction $S^{\oplus,Q}\cc M$ similar
to~\cite[section~6]{GH}. Then the $n$th space of Grayson's
$K$-theory spectrum is given by
   $$K^{Gr}(\cc M)_n=|\Ob S^{\oplus,Q}\cc M|,$$
where $Q=\{1,\ldots,n\}$. It is verified similar
to~\cite[section~6]{GH} that $K^{Gr}(\cc M)$ is a symmetric
spectrum.

If we consider additive categories $\cc P'(X,Y)$, $X,Y\in Sm/F$,
together with associative biexact functors
   $$\cc P'(X,Y)\times\cc P'(Y,Z)\to\cc P'(X,Z),$$
then we shall obtain a spectral category $\cc O_{K^{Gr}}$ on $Sm/F$
such that $\cc O_{K^{Gr}}(X,Y)=K^{Gr}(\cc P'(X,Y))$. The natural map
described above $S^\oplus\cc P'(X,Y)\to S\cc P'(X,Y)$, $X,Y\in
Sm/F$, determines a stable equivalence of symmetric spectra
   $$\cc O_{K^{Gr}}(X,Y)\to\cc O_{K^{\oplus}}(X,Y).$$
Altogether these maps give an equivalence of spectral categories
   $$\cc O_{K^{Gr}}\to\cc O_{K^{\oplus}}.$$

Let $\cc O$ be a spectral category. The category $\Mod\cc O$ of $\cc
O$-modules is enriched over symmetric spectra. Namely, to any
$M,N\in\Mod\cc O$ we associate the symmetric spectrum
   $$Sp^\Sigma(M,N):=\int_{o\in\Ob\cc O}\underline{Sp}^\Sigma(M(o),N(o)),$$
where the integral stands for the coend and
$\underline{Sp}^\Sigma(-,-)$ stands for the internal symmetric
spectrum (see~\cite[section~2.2]{DRO}). By the ``Enriched Yoneda
Lemma" there is a natural isomorphism of symmetric spectra
   $$M(o)\cong Sp^\Sigma(\cc O(-,o),M)$$
for all $o\in\Ob\cc C$ and $M\in\Mod\cc O$.

Let $\cc O$ be symmetric monoidal and let $\diamond:\cc O\wedge\cc
O\to\cc O$ be the structure spectral functor (see
Definition~\ref{avtobus}(4)). By a theorem of Day~\cite{Day}
$\Mod\cc O$ is a closed symmetric monoidal category with smash
product $\wedge$ and $\cc O(-,u)$ being the monoidal unit. The smash
product is defined as
   \begin{equation}\label{smash}
    M\wedge_{\cc O} N=\int^{\Ob\cc O\otimes\cc O}M(o)\wedge N(p)\wedge\cc O(-,o\diamond p).
   \end{equation}
The internal Hom functor, right adjoint to $-\wedge_{\cc O}M$, is
given by
   $$\underline{\Mod}\cc O(M,N)(o):=Sp^\Sigma(M,N(o\diamond-))=\int_{p\in\Ob\cc O}\underline{Sp}^\Sigma(M(p),N(o\diamond p)).$$
It follows from~\cite[2.7]{DRO} that there is a natural isomorphism
   $$\cc O(-,o)\wedge_{\cc O}\cc O(-,p)\cong\cc O(-,o\diamond p).$$

A {\it morphism\/} $\Psi:\cc O\to\cc R$ of spectral categories is
simply a spectral functor. The {\it restriction of scalars\/}
   $$\Psi^*:\Mod\cc R\to\Mod\cc O,\quad M\mapsto M\circ\Psi$$
has a left adjoint functor $\Psi_*$, also denoted $-\wedge_{\cc
O}\cc R$, which we refer to as extension of scalars. It is given by
an enriched coend, i.e., for an $\cc O$-module $N$ the $\cc
R$-module $\Psi_*N=N\wedge_{\cc O}\cc R$ is
   $$\int^{o\in\Ob\cc O}N(o)\wedge\cc R(-,\Psi(o)).$$

Recall that the {\it underlying category\/} $\cc U\cc O$ of a
spectral category $\cc O$ has the same objects as $\cc O$ and the
Hom-sets are defined as $\Hom_{\cc U\cc
O}(o,o')=\Hom_{Sp^\Sigma}(S,\cc O(o,o'))$. Suppose $\cc C$ is a
small category and
   $$f:\cc C\to\cc U\cc O$$
a functor. Denote by $Pre^\Sigma(\cc C)$ the category of presheaves
of symmetric spectra on $\cc C$. Let $U:\Mod\cc O\to Pre^\Sigma(\cc
C)$ be the forgetful functor. It can be proved similar
to~\cite[X.4.1]{Mac} that $U$ has a left adjoint $F:Pre^\Sigma(\cc
C)\to\Mod\cc O$ defined as
   $$F(M)=\int^{c\in\Ob\cc C}M(c)\wedge\cc O(-,f(c)),\quad M\in Pre^\Sigma(\cc C).$$

\section{Model category structures for $\Mod\cc O$}

Let $\cc O$ be a spectral category and let $\Mod\cc O$ be the
category of $\cc O$-modules. We refer the reader to~\cite{Dug} for
basic facts about model categories enriched over symmetric spectra.

\begin{defs}{\rm
(1) A morphism $f$ in $\Mod\cc O$ is a

\begin{itemize}
\item[$\diamond$] {\it level weak equivalence\/} if $f(c)$ is a level weak equivalence in
$Sp^\Sigma$ for all $c\in\Ob\cc O$.

\item[$\diamond$] {\it projective (flat) level fibration\/} if $f(c)$ is a projective (flat) level
fibration in $Sp^\Sigma$ for all $c\in\Ob\cc O$.

\item[$\diamond$] {\it projective (flat) cofibration\/} if $f$ has the left
lifting property with respect to all projective (flat) level acyclic
fibrations.

\end{itemize}

(2) A morphism $f$ in $\Mod\cc O$ is a

\begin{itemize}
\item[$\diamond$] {\it stable weak equivalence\/} if $f(c)$ is a stable weak equivalence in
$Sp^\Sigma$ for all $c\in\Ob\cc O$.

\item[$\diamond$] {\it stable projective (flat) level fibration\/} if $f(c)$ is a stable projective (flat)
fibration in $Sp^\Sigma$ for all $c\in\Ob\cc O$.

\item[$\diamond$] {\it stable projective (flat) cofibration\/} if $f$ has the left
lifting property with respect to all stable projective (flat)
acyclic fibrations.

\end{itemize}

}\end{defs}

\begin{thm}[\cite{DRO,SS}]\label{modelo}
The category $\Mod\cc O$ admits the following four model structures:

\begin{itemize}
\item[$\diamond$] In the projective (respectively flat) level model structure the
weak equivalences are the level weak equivalences and fibrations are
the projective (respectively flat) level fibrations.

\item[$\diamond$] In the projective (respectively flat) stable model structure the
weak equivalences are the stable weak equivalences and fibrations
are the stable projective (respectively flat) fibrations.
\end{itemize}
The four model structures are cellular, proper, spectral and weakly
finitely generated. Moreover, if $\cc O$ is a symmetric monoidal
spectral category then each of the four model structures on $\Mod\cc
O$ is symmetric monoidal with respect to the smash
product~\eqref{smash} of $\cc O$-modules and satisfies the monoid
axiom.
\end{thm}

\begin{proof}
Let us consider one of the four model structures of symmetric
spectra stated in Theorems~\ref{spect}-\ref{spectra}. By those
theorems $Sp^\Sigma$ is a weakly finitely generated monoidal model
category and the monoid axiom holds in $Sp^\Sigma$. It follows
from~\cite[4.2]{DRO} that $\Mod\cc O$ is a cofibrantly generated,
weakly finitely generated model category. Let $I$ and $J$ be the
family of generating cofibrations and trivial cofibrations
respectively. Then the sets of maps in $\Mod\cc O$
   \begin{equation*}\label{ccpi}
    \cc P_I=\{\cc O(-,o)\wedge si\xrightarrow{\cc O(-,o)\wedge i}\cc O(-,o)\wedge ti\mid i\in I,o\in\Ob\cc O\}
   \end{equation*}
and
   $$\cc P_J=\{\cc O(-,o)\wedge sj\xrightarrow{\cc O(-,o)\wedge j}\cc O(-,o)\wedge tj\mid j\in J,o\in\Ob\cc O\}$$
are families of generating cofibrations and trivial cofibrations
respectively. Since limits of $\cc O$-modules are formed
objectwise~\cite[2.2]{DRO} and $Sp^\Sigma$ is cellular by
Theorems~\ref{spect}-\ref{spectra}, then $\Mod\cc O$ is cellular. It
is a spectral model category by~\cite[A.1.1]{SS}
and~\cite[4.4]{DRO}. It is right proper by~\cite[4.8]{DRO}. The
model structure is also left proper because cofibrations are in
particular objectwise monomorphisms, and pushouts along
monomorphisms preserve level/stable weak equivalences of symmetric
spectra.

Finally, if $\cc O$ is a symmetric monoidal spectral category then
by~\cite[4.4]{DRO} $\Mod\cc O$ is symmetric monoidal with respect to
the smash product~\eqref{smash} of $\cc O$-modules and satisfies the
monoid axiom.
\end{proof}

\section{Motivic model category structures}

In what follows, if otherwise is specified, we work with spectral
categories $\cc O$ over $Sm/F$ such that there is a functor of
categories
   \begin{equation}\label{funct}
    u:Sm/F\to\cc U\cc O,
   \end{equation}
which is identical on objects. One has a bifunctor
   $$\cc O(-,-):Sm/F^{\op}\times Sm/F\to Sp^\Sigma.$$
$\cc O_{K^{\oplus}},\cc O_{K},\cc O_{K^{Gr}},\cc
O_{K^{\oplus}_0},\cc O_{K_0},\cc O_{cor}$ are examples of such
spectral categories. Note that $\cc O$ is a $\cc O_{naive}$-algebra
in the sense that there is a spectral functor $\cc O_{naive}\to\cc
O$ induced by the functor $u$. The spectral category $\cc O_{naive}$
plays the same role as the ring of integers for abelian groups or
the sphere spectrum for symmetric spectra.

Regarding $\cc O_{naive}$-modules as presheaves (see
Example~\ref{urr}) of symmetric spectra $Pre^\Sigma(Sm/F)$, we get a
pair of adjoint functors
   \begin{equation}\label{polezno}
    \xymatrix{{\Psi_*}:Pre^\Sigma(Sm/F)\ar@<0.5ex>[r]&\Mod\cc O:{\Psi^*}.\ar@<0.5ex>[l]}
   \end{equation}
One has for all $X\in Sm/F$,
   \begin{equation}\label{polez}
    \Psi_*(\cc O_{naive}(-,X))=\cc O_{naive}(-,X)\wedge_{\cc O_{naive}}\cc O\cong\cc O(-,X).
   \end{equation}

For simplicity we work with the stable projective model structure on
$\Mod\cc O$ from now on. The interested reader can also consider the
stable flat model structure as well.

Recall that the Nisnevich topology is generated by elementary
distinguished squares, i.e. pullback squares
   \begin{equation}\label{squareQ}
    \xymatrix{\ar@{}[dr] |{\textrm{$Q$}}U'\ar[r]\ar[d]&X'\ar[d]^\phi\\
              U\ar[r]_\psi&X}
   \end{equation}
where $\phi$ is etale, $\psi$ is an open embedding and
$\phi^{-1}(X\setminus U)\to(X\setminus U)$ is an isomorphism of
schemes (with the reduced structure). Let $\cc Q$ denote the set of
elementary distinguished squares in $Sm/F$. By $\cc Q_{\cc O}$
denote the set of squares
   \begin{equation}\label{squareOQ}
    \xymatrix{\ar@{}[dr] |{\textrm{$\cc O Q$}}\cc O(-,U')\ar[r]\ar[d]&\cc O(-,X')\ar[d]^\phi\\
              \cc O(-,U)\ar[r]_\psi&\cc O(-,X)}
   \end{equation}
which are obtained from the squares in $\cc Q$ by taking $X\in Sm/F$
to $\cc O(-,X)$. The arrow $\cc O(-,U')\to\cc O(-,X')$ can be
factored as a cofibration $\cc O(-,U')\rightarrowtail Cyl$ followed
by a simplicial homotopy equivalence $Cyl\to\cc O(-,X')$. There is a
canonical morphism $A_{\cc O Q}:=\cc O(-,U)\coprod_{\cc O(-,U')}
Cyl\to\cc O(-,X)$.

\begin{defs}{\rm
(1) The {\it Nisnevich local model structure\/} on $\Mod\cc O$ is
the Bousfield localization of the stable projective model structure
(see Theorem~\ref{modelo}) with respect to the set of projective
cofibrations
   \begin{equation*}\label{no}
    \cc N_{\cc O}=\{\cyl(A_{\cc O Q}\to\cc O(-,X))\}_{\cc Q_{\cc O}}.
   \end{equation*}
The homotopy category for the Nisnevich local model structure will
be denoted by $\shnis\cc O$. If $\cc O=\cc O_{naive}$ then we shall
write $\shnis(F)$ to denote $\shnis\cc O_{naive}$.

(2) The {\it motivic model structure\/} on $\Mod\cc O$ is the
Bousfield localization of the Nisnevich local model structure with
respect to the set of projective cofibrations
   \begin{equation*}\label{ao}
    \cc A_{\cc O}=\{\cyl(\cc O(-,X\times\bb A^1)\to\cc O(-,X))\}_{X\in Sm/F}.
   \end{equation*}
The homotopy category for the motivic model structure will be
denoted by $\sheff\cc O$. If $\cc O=\cc O_{naive}$ then we shall
write $\sheff(F)$ to denote $\sheff\cc O_{naive}$.

}\end{defs}

\begin{rem}{\rm
A stably fibrant $\cc O$-module $M$ is Nisnevich local \ifff it is
{\it flasque\/} in the sense that for each elementary distinguished
square $Q\in\cc Q$ the square of symmetric spectra $M(Q)$ is
homotopy pullback.

}\end{rem}

Before collecting properties for the Nisnevich local and motivic
model structures we need to recall some facts from unstable $\bb
A^1$-topology.

Let $Pre^{\bb N}(Sm/F)$ be the category of presheaves of ordinary
simplicial spectra $Sp^{\bb N}$. Then $Sp^{\bb N}$ and $Pre^{\bb
N}(Sm/F)$ enjoy the stable projective model structure defined
similar to $Sp^\Sigma$ and $Pre^{\Sigma}(Sm/F)$ (see~\cite{Hov}).
By~\cite[3.5]{Hov} the stable projective model structure on $Sp^{\bb
N}$ coincides with the stable model structure of
Bousfield--Friedlander~\cite{BF}. By~\cite[4.2.5]{HSS} the forgetful
functor $U:Sp^\Sigma\to Sp^{\bb N}$ has a left adjoint $V$ and the
pair $(U,V)$ forms a Quillen equivalence of the stable model
categories. The Nisnevich local and motivic model structures on
$Pre^{\bb N}(Sm/F)$ are defined similar to those on
$Pre^\Sigma(Sm/F)$ by means of the Bousfield localization.

\begin{lem}\label{modelpres}
The adjoint pair $(V,U):Sp^{\bb N}\leftrightarrows Sp^\Sigma$ can be
extended to an adjoint pair $(V,U):Pre^{\bb N}(Sm/F)\leftrightarrows
Pre^\Sigma(Sm/F)$. It forms a Quillen equivalence for the stable
projective, Nisnevich local and motivic model structures
respectively.
\end{lem}

\begin{proof}
See~\cite[section~10]{Hov} and~\cite[section~4.5]{Jar2}.
\end{proof}

Let $Pre(Sm/F)$ denote the category of pointed simplicial presheaves
on $Sm/F$. Then $Pre(Sm/F)$ enjoys the projective model structure in
which fibrations and weak equivalences are defined schemewise.
Projective cofibrations are those maps which have the corresponding
lifting property. As above, one defines the Nisnevich local
projective and motivic model structures by means of the Bousfield
localization. By~\cite{Bl} both model structures are proper
simplicial cellular and the Nisnevich local projective model
structure coincides with the model structure in which weak
equivalences are the local weak equivalences with respect to
Nisnevich topology and cofibrations are the projective cofibrations.

In~\cite{Hov} Hovey constructed stable (symmetric) model structures
out of certain model categories with certain Quillen endofunctors.
For example, one can apply Hovey's constructions to Nisnevich local
and motivic projective model structures on $Pre(Sm/F)$ with $-\wedge
S^1$ a Quillen endofunctor. Denote the resulting model categories by
$Sp^{\bb N}_{nis}(Pre(Sm/F))$ and $Sp^{\bb N}_{mot}(Pre(Sm/F))$
(respectively $Sp^{\Sigma}_{nis}(Pre(Sm/F))$ and
$Sp^{\Sigma}_{mot}(Pre(Sm/F))$). As categories $Sp^{\bb
N}_{nis}(Pre(Sm/F))$ and $Sp^{\bb N}_{mot}(Pre(Sm/F))$ (respectively
$Sp^{\Sigma}_{nis}(Pre(Sm/F))$ and $Sp^{\Sigma}_{mot}(Pre(Sm/F))$)
coincide with $Pre^{\bb N}(Sm/F)$ (respectively
$Pre^{\Sigma}(Sm/F)$). The following proposition states that the
corresponding model category structures on these coincide as well.

\begin{prop}\label{modelcomp}
The model categories $Sp^{\bb N}_{nis}(Pre(Sm/F))$ and $Sp^{\bb
N}_{mot}(Pre(Sm/F))$ (respectively $Sp^{\Sigma}_{nis}(Pre(Sm/F))$
and $Sp^{\Sigma}_{mot}(Pre(Sm/F))$) coincide with the Nisnevich
local projective and motivic model category structures on $Pre^{\bb
N}(Sm/F)$ ($Pre^{\Sigma}(Sm/F)$) respectively.
\end{prop}

\begin{proof}
We prove the statement for the Nisnevich local projective model
structure. The statement for the motivic model structure is proved
in a similar way.

By Hovey's construction of the model category $Sp^{\bb
N}_{nis}(Pre(Sm/F))$ a map of spectra $i:A\to B$ is a cofibration
\ifff the induced maps $i_0:A_0\to B_0$ and
$j_n:A_n\coprod_{A_{n-1}\wedge S^1}B_{n-1}\wedge S^1\to B_n$, $n\geq
1$, are projective cofibrations in $Pre(Sm/F)$. It follows that
cofibrations in $Sp^{\bb N}_{nis}(Pre(Sm/F))$ coincide with
cofibrations in the stable projective model structure on $Pre^{\bb
N}(Sm/F)$, and hence with cofibrations in the Nisnevich local model
structure on $Pre^{\bb N}(Sm/F)$ because the Bousfield localization
preserves cofibrations.

An object $X$ in $Sp^{\bb N}_{nis}(Pre(Sm/F))$ is fibrant \ifff it
is levelwise fibrant, i.e. levelwise projective fibrant and flasque,
and each map $X_n\to\Omega X_{n+1}$, $n\geq 0$, which is adjoint to
the structure map $X_n\wedge S^1\to X_{n+1}$, is a weak equivalence.
It follows that each map $X_n\to\Omega X_{n+1}$, $n\geq 0$, is a
schemewise weak equivalence, and hence $X$ is stable fibrant in the
stable projective model structure on $Pre^{\bb N}(Sm/F)$. Since $X$
is levelwise flasque, we see that $X(Q)$ is homotopy cartesian for
every elementary distinguished square $Q\in\cc Q$. Therefore $X$ is
fibrant in the Nisnevich local model structure on $Pre^{\bb
N}(Sm/F)$. One easily sees that every fibrant object in the
Nisnevich local model structure on $Pre^{\bb N}(Sm/F)$ is fibrant in
$Sp^{\bb N}_{nis}(Pre(Sm/F))$. We have shown that cofibrations and
fibrant objects in both model categories coincide.

By~\cite[9.7.4(4)]{Hir} a map $g$ in a simplicial model category is
a weak equivalence if and only if for some cofibrant approximation
$\tilde g:\wt X\to\wt Y$ to $g$ and every fibrant object $Z$ the map
of simplicial sets $\tilde g^*:\Map(Y,Z)\to\Map(X,Z)$ is a weak
equivalence. We infer that weak equivalences in $Sp^{\bb
N}_{nis}(Pre(Sm/F))$ and in the Nisnevich local model structure on
$Pre^{\bb N}(Sm/F)$ coincide. Thus the model structure on $Sp^{\bb
N}_{nis}(Pre(Sm/F))$ coincides with the Nisnevich local model
structure on $Pre^{\bb N}(Sm/F)$. The fact that
$Sp^{\Sigma}_{nis}(Pre(Sm/F))$ and the Nisnevich local model
structure on $Pre^{\Sigma}(Sm/F)$ coincide is checked similar to
presheaves of ordinary spectra.
\end{proof}

In order to show some homotopically important properties for
$Pre^{\Sigma}(Sm/F)$, we need to discuss Jardine's model structures.
$Pre(Sm/F)$ enjoys the injective model structure in which
cofibrations are monomorphisms and weak equivalences are the local
weak equivalences. Global fibrations are those maps which have the
right lifting property with respect to all maps which are
cofibrations and local weak equivalences. As above, one defines the
motivic injective model structure by means of the Bousfield
localization. Then both the Nisnevich local injective and motivic
injective model structures satisfy the axioms for a proper
simplicial model category~\cite{Jar1,Jar2}. One easily sees that the
Nisnevich local projective (respectively motivic projective) model
structure on $Pre(Sm/F)$ is Quillen equivalent to the Nisnevich
local injective (respectively motivic injective) model structure.

Denote by $Sp^\Sigma_{nis,J}(Pre(Sm/F))$ and
$Sp^\Sigma_{mot,J}(Pre(Sm/F))$ (``$J$" for Jardine) stable symmetric
model structures corresponding to the Nisnevich local injective
(respectively motivic injective) model structure on $Pre(Sm/F)$.
By~\cite[4.32]{Jar2} and~\cite[Thm.~12]{Jar3} these are proper
simplicial model categories. It follows from~\cite[9.3]{Hov} that
the natural functors
   $$Sp^\Sigma_{nis}(Pre(Sm/F))\to Sp^\Sigma_{nis,J}(Pre(Sm/F))$$
and
   $$Sp^\Sigma_{mot}(Pre(Sm/F))\to Sp^\Sigma_{mot,J}(Pre(Sm/F))$$
are Quillen equivalences.

The category $Pre^\Sigma(Sm/F)$ of presheaves of symmetric spectra
also enjoys the following model structures (see~\cite{Jar2,Jar3} for
details). A map $f:X\to Y$ is a {\it level equivalence\/} if each of
the component maps $f:X_n\to Y_n$ is a local weak equivalence
(respectively motivic equivalence). The map $f$ is a {\it level
cofibration\/} if each of the maps $X_n\to Y_n$ is a monomorphism of
simplicial presheaves. Denote these model categories by
$Pre^\Sigma_{nis}(Sm/F)$ and $Pre_{mot}^\Sigma(Sm/F)$ respectively.
These are proper simplicial model categories.

For every $X\in Pre^\Sigma(Sm/F)$ we obtain a natural construction
   $$X\lra{i_1}X_s\lra{i_2}X_{si}$$
of an {\it injective stably fibrant model\/} $X_{si}$, where $i_1$
is a trivial stable cofibration in the model category
$Sp^\Sigma_{nis,J}(Pre(Sm/F))$ (respectively in
$Sp^\Sigma_{mot,J}(Pre(Sm/F))$) and $i_2$ is a level cofibration and
a level equivalence in $Pre_{nis}^\Sigma(Sm/F)$ (respectively in
$Pre_{mot}^\Sigma(Sm/F)$).

By~\cite{Jar2,Jar3} a map $f:X\to Y$ of $Pre^\Sigma(Sm/F)$ is a
stable weak equivalence \ifff it induces a weak equivalence of
simplicial mapping Kan complexes
   $$f:\Map(Y,W)\to\Map(X,W)$$
for each injective stably fibrant object $W$. Observe that the maps
$i_1$ and $i_2$ are both stable weak equivalences.

\begin{prop} \textrm{\em (\cite[4.41]{Jar2})}\label{jardine}
Suppose that $i:A\to B$ is a stable cofibration in
$Sp^\Sigma_{nis,J}(Pre(Sm/F))$ or in $Sp^\Sigma_{mot,J}(Pre(Sm/F))$
and that $j:C\to D$ is a level cofibration. Then the map
   $$(i,j)_*:(B\wedge_{\cc O_{naive}}C)\cup_{(A\wedge_{\cc O_{naive}}C)}(A\wedge_{\cc O_{naive}}D)
     \to B\wedge_{\cc O_{naive}}D$$
is a level cofibration. If either $i$ or $j$ is a stable equivalence
in $Sp^\Sigma_{nis,J}(Pre(Sm/F))$ (respectively in
$Sp^\Sigma_{mot,J}(Pre(Sm/F))$), then so is $(i,j)_*$.
\end{prop}

The proposition was actually shown for
$Sp^\Sigma_{mot,J}(Pre(Sm/F))$. However the proof of this result is
really quite generic, and holds essentially anywhere that one
succeeds in generating the usual machinery of symmetric spectrum.
This includes the present discussion of symmetric $S^1$-spectra in
the Nisnevich local case $Sp^\Sigma_{nis,J}(Pre(Sm/F))$.

\begin{thm}[Jardine]\label{uhh}
The Nisnevich local projective and motivic projective model
structures on the category $Pre^\Sigma(Sm/F)$ of presheaves of
symmetric spectra are cellular, proper, spectral, weakly finitely
generated, symmetric monoidal and satisfy the monoid axiom.
\end{thm}

\begin{proof}
Since Bousfield localization respects cellularity and left
properness, then both model structures are cellular and left proper.
Right properness is proved similar to~\cite[4.15]{Jar2}
and~\cite[Thm.~12]{Jar3}.

It follows from~\cite[3.2.13]{HSS} that every object of
$Pre^\Sigma(Sm/F)$ is small. By Lemma~\ref{pqp} both model
structures are weakly finitely generated, because domains and
codomains of morphisms from $\cc N\cup\cc A$ are finitely
presentable. These are plainly spectral as well.
Proposition~\ref{jardine} and Theorem~\ref{modelo} imply that both
model structures are symmetric monoidal.

It remains to verify the monoid axiom. Let $i:A\to B$ be a trivial
stable cofibration and let $C\in Pre^\Sigma(Sm/F)$. By
Proposition~\ref{jardine} the map
   $$i\wedge 1:A\wedge_{\cc O_{naive}}C\to B\wedge_{\cc O_{naive}}C$$
is a level cofibration and a stable equivalence. Therefore for any
injective stably fibrant object $W$ the map $\Map(i\wedge 1,W)$ is a
trivial fibration of simplicial sets. Let $j$ be a pushout of
$i\wedge 1$ along some map. It follows that $\Map(j,W)$ is a trivial
fibration of simplicial sets, and hence $j$ is a stable equivalence.
Now the monoid axiom follows from the fact that any transfinite
composition of weak equivalences is a weak equivalence in a weakly
finitely generated model category.
\end{proof}

Given a spectral category $\cc O$ over $Sm/F$, we want to establish
the same properties for the Nisnevich local and motivic model
structures on $\Mod\cc O$ as for $Pre^\Sigma(Sm/F)$ from the
preceding theorem. For this we have to give the following

\begin{defs}{\rm
(1) We say that $\cc O$ is {\it Nisnevich excisive\/} if for every
elementary distinguished square $Q$
    \begin{equation*}
    \xymatrix{\ar@{}[dr] |{\textrm{$Q$}}U'\ar[r]\ar[d]&X'\ar[d]^\phi\\
              U\ar[r]_\psi&X}
   \end{equation*}
the square $\cc O Q$~\eqref{squareOQ} is homotopy pushout in the
Nisnevich local model structure on $Pre^\Sigma(Sm/F)$.

(2) $\cc O$ is {\it motivically excisive\/} if:

\begin{itemize}
\item[(A)] for every elementary distinguished square $Q$ the square $\cc O
Q$~\eqref{squareOQ} is homotopy pushout in the motivic model
structure on $Pre^\Sigma(Sm/F)$ and

\item[(B)] for every $X\in Sm/F$ the natural map
   $$\cc O(-,X\times\bb A^1)\to\cc O(-,X)$$
is a weak equivalence in the motivic model structure on
$Pre^\Sigma(Sm/F)$.
\end{itemize}

}\end{defs}

The following lemma says that property~(B) is redundant for
symmetric monoidal spectral categories.

\begin{lem}\label{qoqo}
Let $\cc O$ be a symmetric monoidal spectral category on $Sm/F$ such
that the monoidal product is given by cartesian product of schemes.
Then the map
   $$f:\cc O(-,X\times\bb A^1)\to\cc O(-,X)$$
is a weak equivalence in the motivic model structure on
$Pre^\Sigma(Sm/F)$.
\end{lem}

\begin{proof}
We follow an argument of~\cite[p.~694]{RO}. As in classical
algebraic topology, an inclusion of motivic spaces $g:A\to B$ is an
$\bb A^1$-deformation retract if there exist a map $r:B\to A$ such
that $rg=\id_A$ and an $\bb A^1$-homotopy $H:B\times\bb A^1\to B$
between $gr$ and $\id_B$ which is constant on $A$. Then $\bb
A^1$-deformation retracts are motivic weak equivalences.

There is an obvious map $r:\cc O(-,X)\to\cc O(-,X\times\bb A^1)$
such that $fr=1$. Since $\cc O$ is a symmetric monoidal spectral
category, it follows that
   $$\cc O(-\times\bb A^1,X\times\bb A^1)\cong\underline\Mod\cc O(\cc O(-,\bb A^1),\cc O(-,X\times\bb A^1))$$
is an $\cc O$-module.

There is a natural isomorphism of symmetric spectra
   $$Sp^\Sigma(\cc O(-,X\times\bb A^1),\cc O(-\times\bb A^1,X\times\bb A^1))
     \cong\cc O(X\times\bb A^1\times\bb A^1,X\times\bb A^1).$$
Consider the functor~\eqref{funct} of categories $u:Sm/F\to\cc U\cc
O$. Denote by $\alpha$ the obvious map $\bb A^1\times\bb A^1\to\bb
A^1$. We set $h=u(1_X\times\alpha)$; then $h$ uniquely determines a
morphism of $\cc O$-modules
   $$h':\cc O(-,X\times\bb A^1)\to\cc O(-\times\bb A^1,X\times\bb A^1).$$
This morphism can be regarded as a morphism of $\cc
O_{naive}$-modules, denoted by the same letter. By adjointness $h'$
uniquely determines a map of $\cc O_{naive}$-modules
   \begin{gather*}
    H:\cc O(-,X\times\bb A^1)\wedge_{\cc O_{naive}}\cc O_{naive}(-,\bb A^1)\to\cc O(-,X\times\bb A^1).
   \end{gather*}
Then $H$ yields a level $\bb A^1$-homotopy between the identity map
and $rf$. We see that $f$ is a level motivic equivalence, and hence
it is a weak equivalence in $Sp^{\Sigma}_{mot}(Pre(Sm/F))$.
\end{proof}

\begin{thm}\label{panin}
$\cc O_{K^{\oplus}},\cc O_{K},\cc O_{K^{Gr}},\cc
O_{K^{\oplus}_0},\cc O_{K_0},\cc O_{cor}$ are Nisnevich excisive
spectral categories.
\end{thm}

\begin{proof}
We want to prove the statement first for $\cc O_{K}$. The cases $\cc
O_{K^\oplus},\cc O_{K^{Gr}}$ are checked in a similar way. Given
$q\in\bb Z$ and a smooth scheme $X$, let $\cc K_q(-,X)$ denote the
sheaf associated to the presheaf $W\mapsto K_q(W,X)=\pi_q(\cc
O_K(W,X))$. If we show that $\cc O_K(Q)$ is homotopy pushout in
$Sp_{\nis}^{\bb N}(Pre(Sm/F))$ for every elementary distinguished
square $Q$ then it will follow from~\cite[Lemma~10]{Jar3} that $\cc
O_K$ is a Nisnevich excisive spectral category.

We have to verify that for any elementary distinguished square $Q$
the sequence of sheaves
   $$\cdots\to\cc K_{q+1}(-,X)\to\cc K_q(-,U')\to\cc K_q(-,U)\ps\cc K_q(-,X')\to\cc K_q(-,X)\to\cdots$$
is exact. Because the Nisnevich topology has enough points, the
sequence
   $$\cdots\to K_{q+1}(-,X)\to K_q(-,U')\to K_q(-,U)\ps K_q(-,X')\to K_q(-,X)\to\cdots$$
will become exact after sheafifying precisely if it becomes exact
whenever one applies the presheaves to the Henselization $W$ of a
smooth scheme $T$ at a point $t$. Thus it is enough to show that for
any such $W$ the square
   \begin{equation}\label{hensel}
    \xymatrix{\cc O_K(W,U')\ar[d]\ar[r]&\cc O_K(W,X')\ar[d]\\
               \cc O_K(W,U)\ar[r]&\cc O_K(W,X)}
   \end{equation}
is homotopy pushout (=homotopy pullback) of spectra. We shall
actually show that~\eqref{hensel} gives a split exact sequence
   \begin{equation}\label{ter}
    \cc O_K(W,U')\hookrightarrow\cc O_K(W,U)\ps\cc O_K(W,X')\twoheadrightarrow\cc O_K(W,U')
   \end{equation}
in the (triangulated) homotopy category of spectra $\Ho Sp$.

Given a triangulated category $\cc T$, consider a commutative
diagram
   \begin{equation}\label{ivan}
    \xymatrix{A'_1\ar[d]_{\pi_1}\ar@/^30pt/[rrr]^{\alpha'}\ar@/^10pt/[rr]^{\alpha'\times 0}\ar[r]_{i'}
             &C'\ar[d]_\delta\ar[r]_(.4){\Phi'\times\Psi'}&A'\times B'\ar[d]^{\pi\times\rho}\ar[r]_(.6){p'}&A'\ar[d]^\pi\\
             A_1\ar@/_30pt/[rrr]^\alpha\ar@/_10pt/[rr]_{\alpha\times 0}\ar[r]^i&C\ar[r]^(.4){\Phi\times\Psi}
             &A\times B\ar[r]^(.6)p&A,}
   \end{equation}
where $\rho:B'\to B$, $\Phi\times\Psi$, $\Phi'\times\Psi'$ are
isomorphisms, $\alpha=\Phi i$ and $\alpha'=\Phi' i'$ are
isomorphisms. Clearly, the right and the middle squares are
cartesian, as well as so is the outer square with vertices
$(A_1',A',A_1,A)$. Therefore the left square is cartesian and,
moreover, the induced sequence
   $$A'_1\to A_1\ps C'\to C$$
is split exact in $\cc T$. Below we shall be constructing such a
diagram in the triangulated category $\Ho Sp$.

Consider an elementary distinguished square $Q$
    \begin{equation*}
    \xymatrix{\ar@{}[dr] |{\textrm{$Q$}}U'\ar[r]^{i'}\ar[d]_{\pi_1}&X'\ar[d]^\pi\\
              U\ar[r]_i&X.}
   \end{equation*}
Denote by $\cc P(W,X)^U$ (respectively $\cc P(W,X)^{\neg U}$) the
full subcategory of $\cc P(W,X)$ consisting of those bimodules
$P\in\cc P(W,X)$ for which $\supp(P)\subseteq W\times U$
(respectively $\supp(P)\varsubsetneq W\times U$). There exists a
functor
   $$\Phi:\cc P(W,X)\to\cc P(W,X)^U,$$
constructed as follows. Given $P\in\cc P(W,X)$ set $S:=\supp(P)$.
Then $S=S'\sqcup S''$ with $S'\subseteq W\times U$ and $S''$ such
that its each connected component is not contained in $W\times U$.
The sheaf $P$ is canonically equal to $P'\ps P''$ with
$\supp(P')=S'$ and $\supp(P'')=S''$. Moreover, if $P_1\in\cc P(W,X)$
and $P_1=P_1'\ps P_1''$ is a similar decomposition, then every
morphism $f:P\to P_1$ is of the form $f=\left(\begin{array}{cc} f' &
0
\\ 0 & f'' \\ \end{array}\right)$. This is because
$\Hom(P',P_1'')=\Hom(P'',P_1')=0$ and the corresponding supports are
disjoint. So we set $\Phi(P)=P',\Phi(f)=f'$. There is also a functor
   $$\Psi:\cc P(W,X)\to\cc P(W,X)^{\neg U},$$
defined as $\Psi(P)=P'',\Psi(f)=f''$. The full subcategories $\cc
P(W,X')^{U'},\cc P(W,X')^{\neg U'}\subset\cc P(W,X')$ and functors
$\Phi':\cc P(W,X')\to\cc P(W,X')^{U'},\Psi':\cc P(W,X')\to\cc
P(W,X')^{\neg U'}$ are defined in a similar way.

The map $\pi:X'\to X$ induces two functors
   $$\pi_*^U:\cc P(W,X')^{U'}\to\cc P(W,X)^U,\quad\pi_*^{\neg U}:\cc P(W,X')^{\neg U'}\to\cc P(W,X)^{\neg U}.$$
Indeed, if $Q\in\cc P(W,X')^{U'}$ then $\supp(Q)\subset W\times U'$.
It follows that $\supp(\pi_*(Q))=\pi(\supp(Q))\subset W\times U$,
and hence $\pi_*(Q)\in\cc P(W,X)^{U}$. One sets
$\pi_*^U:=\pi_*|_{\cc P(W,X')^{U'}}$. If $Q\in\cc P(W,X')^{\neg U'}$
then no connected component of $\supp(Q)$ is contained in $W\times
U'$. Since $\pi^{-1}(U)=U'$, it follows that no connected component
of $\supp(\pi_*(Q))=\pi(\supp(Q))$ is contained in $W\times U$. We
see that $\pi_*(Q)\in\cc P(W,X)^{\neg U}$ and one puts $\pi_*^{\neg
U}:=\pi_*|_{\cc P(W,X')^{\neg U'}}$.

Consider a diagram of categories
   \begin{equation}\label{sup}
    \xymatrix{\cc P(W,U')\ar[d]_{\pi_1,*}\ar[r]^{i'_*}&\cc P(W,X')\ar[d]_{\pi_*}\ar[r]^(.33){\Phi'\times\Psi'}
    &\cc P(W,X')^{U'}\times\cc P(W,X')^{\neg U'}\ar[d]^{\pi_*^U\times\pi_*^{\neg U}}\ar[r]^(.65){p'}
    &\cc P(W,X')^{U'}\ar[d]^{\pi_*^U}\\
    \cc P(W,U)\ar[r]_{i_*}&\cc P(W,X)\ar[r]_(.33){\Phi\times\Psi}
    &\cc P(W,X)^{U}\times\cc P(W,X)^{\neg U}\ar[r]_(.65)p&\cc P(W,X)^U,}
   \end{equation}
We claim that:
\begin{enumerate}
\item $\Phi\times\Psi$ and $\Phi'\times\Psi'$ are equivalences of
categories;
\item $\Psi\circ i_*$ is the zero functor;
\item $\Psi'\circ i'_*$ is the zero functor;
\item $\Phi'\circ i_*'$ is an equivalence of
categories;
\item $\Phi\circ i_*$ is an equivalence of
categories;
\item $\pi_*^{\neg U}$ is an equivalence of
categories;
\item the diagram is commutative up to natural isomorphisms of
functors.
\end{enumerate}
Statements (1)-(7) together with~\cite[1.3.1]{Wal} will imply that
the diagram of spectra
   $$\xymatrix{\cc O_K(W,U')\ar[d]_{\pi_1,*}\ar[r]^{i'_*}&\cc O_K(W,X')\ar[d]_{\pi_*}\ar[r]^(.33){\Phi'\times\Psi'}
             &\cc O_K(W,X')^{U'}\times\cc O_K(W,X')^{\neg U'}\ar[d]^{\pi_*^U\times\pi_*^{\neg U}}\ar[r]^(.65){p'}
             &\cc O_K(W,X')^{U'}\ar[d]^{\pi_*^U}\\
             \cc O_K(W,U)\ar[r]_{i_*}&\cc O_K(W,X)\ar[r]_(.33){\Phi\times\Psi}
             &\cc O_K(W,X)^{U}\times\cc O_K(W,X)^{\neg U}\ar[r]_(.65)p&\cc O_K(W,X)^U,}$$
obtained from~\eqref{sup} by taking realizations, is commutative in
$\Ho Sp$. So we shall obtain a diagram of the form~\eqref{ivan},
which will yield a split exact sequence~\eqref{ter}, and hence
square~\eqref{hensel} is homotopy pushout, as required.

So it remains to show (1)-(7). Let us show that $\Phi\times\Psi$ is
an equivalence of categories (the same fact for $\Phi'\times\Psi'$
is checked in a similar way). Consider a functor
   $$\Theta:\cc P(W,X)^{U}\times\cc P(W,X)^{\neg U}\to\cc P(W,X),$$
defined as $\Theta(P',P'')=P'\ps P''$, $\Theta(f',f'')=f'\ps f''$.
Clearly, the canonical morphism
$can_P:P\to\Theta\circ(\Phi\times\Psi)(P)$ is an isomorphism for
every $P\in\cc P(W,X)$
   $$\xymatrix{P\ar@/_/[dr]^(.5){\cong}_{can_P}\ar[r]^(.35){\Phi\times\Psi}&(P',P'')\ar[d]^\Theta\\
               &P'\ps P''.}$$
Given a morphism $f:P\to P_1$ in $\cc P(W,X)$, we have that
$\Hom(P'',P_1')=\Hom(P',P_1'')=0$ and the diagram
   $$\xymatrix{P\ar[d]_f\ar[r]^(.35){can_P}&P'\ps P''\ar[d]^{f'\ps f''}\\
               P_1\ar[r]_(.4){can_{P_1}}&P_1'\ps P_1''}$$
is commutative. The latter shows that there is a natural
transformation of functors $can:\id\to\Theta\circ(\Phi\times\Psi)$.
Since $can_P$ is an isomorphism for all $P$, then $can$ is an
isomorphism of functors.

The composition $(\Phi\times\Psi)\circ\Theta$ is just the identity
functor. Indeed,
$[(\Phi,\Psi)\circ\Theta](P',P'')=(\Phi\times\Psi)(P'\ps
P'')=(P',P'')$. So (1) is verified.

For any $P\in\cc P(W,U)$ one has $\supp(i_*(P))=i(\supp(P))\subset
W\times U$. We see that $\Psi(i_*(P))=0$. So (2) is verified.
Property~(3) is checked in a similar way.

Let us prove (5). Consider the functor
   $$i^*|_{\cc P(W,X)^U}:\cc P(W,X)^U\to\cc P(W,U).$$
First, it is well defined. Indeed,
$\supp(i^*(P))=i^{-1}(\supp(P))=\supp(P)$, because $\supp(P)\subset
W\times U$ for all $P\in\cc P(W,X)^U$. For brevity we shall write
$i^*$ instead of $i^*|_{\cc P(W,X)^U}$. Second, one has adjunction
morphisms
   $$adj_P:i^*i_*(P)\to P,\quad P\in\cc P(W,U),$$
and
   $$adj_{P'}:P'\to i_*i^*(P'),\quad P'\in\cc P(W,X)^U.$$
Clearly, these determine natural transformations of functors
$i^*i_*\to\id$ and $\id\to i_*i^*$.

We want to show that these are isomorphisms. If $P\in\cc P(W,U)$
then $\supp(i^*i_*(P))=i^{-1}i(\supp(P))=\supp(P)$. Moreover, for
every $x\in\supp(P)$ the induced morphism of stalks
   $$(i^*i_*(P))_x\to P_x$$
is an isomorphism. Therefore $adj_P$ is an isomorphism. If $P'\in\cc
P(W,X)^U$ it follows that
$\supp(i_*i^*(P'))=i(i^{-1}(\supp(P')))=i(\supp(P'))=\supp(P')$.
Moreover, $adj_{P'}$ induces an isomorphism of stalks
   $$P_x'\to(i_*i^*(P'))_x.$$
Therefore $adj_{P'}$ is an isomorphism. So (5) is verified. Property
(4) is verified in a similar way.

Next, let us check (7). If $Q\in\cc P(W,X')$ and $Q=Q'\ps Q''$ is
its canonical decomposition with $Q'=\Phi'(Q),Q''=\Psi'(Q)$. Then
$\pi_*(Q)=\pi_*(Q')\ps\pi_*(Q'')$. On the other hand,
$\pi_*(Q)=\pi_*(Q)'\ps\pi_*(Q)''$. Comparing supports, one gets that
$\pi_*^U(Q'):=\pi_*(Q')=\pi_*(Q)'$ and $\pi_*^{\neg
U}(Q''):=\pi_*(Q'')=\pi_*(Q)''$. So,
   $$(\pi_*(Q)',\pi_*(Q)'')=(\pi_*^U(Q'),\pi_*^{\neg U}(Q'')).$$
We see that $\pi_*^U\circ\Phi'=\Phi\circ\pi_*$ and $\pi_*^{\neg
U}\circ\Psi'=\Psi\circ\pi_*$. So (7) is verified.

To prove (6), we shall need some notation. Let $Q$ be a coherent
$\cc O_{W\times X'}$-module such that $S(Q):=\supp(Q)$ is
quasi-finite over $W$. Let $S(Q)''\subset S(Q)$ denote the union of
those connected components in $S(Q)$, each of which is not contained
in $W\times U'$. Let $S(Q)'\subset S(Q)$ denote the union of the
other connected components in $S(Q)$. Clearly, $S(Q)'\subset W\times
U'$ and $S(Q)=S(Q)'\sqcup S(Q)''$. Then one has a decomposition
   $$Q=Q'\ps Q''$$
such that $\supp(Q')=S(Q)'$ and $\supp(Q'')=S(Q)''$.

Given a coherent $\mathcal O_{W \times X}$-module $\mathcal M$ which
is coherent as an $\mathcal O_{W}$-module,
we set $\pi^*_{\neg U}(\mathcal M):=\pi^*(\mathcal M)''$.

\begin{sublem}
The following statements are true:

$(a)$ If $P\in\cc P(W,X)^{\neg U}$ then $\pi^*_{\neg U}(P)\in\cc
P(W,X^{\prime})^{\neg U^{\prime}}$. In particular,
$\supp(\pi_*^{\neg U}(P))$ is finite and surjective over $W$.

$(b)$ The composition
$P\xrightarrow{adj}\pi^*\pi_*(P)\to\pi_*(\pi^*_{\neg U}(P))$ is an
isomorphism for all $P\in\cc P(W,X)^{\neg U}$.

$(c)$ The composition $\pi^*_{\neg
U}(\pi_*(Q))\to\pi^*\pi_*(Q)\xrightarrow{adj} Q$ is an isomorphism
for all $Q\in\cc P(W,X')^{\neg U'}$.

\end{sublem}

\begin{proof}
Firstly prove (a). It is easy to check that $\pi_*^{\neg U}(P)$ is
coherent as an $\mathcal O_W$-module. Since $\pi$ is \'{e}tale, it
is flat as well. Thus $\pi_*^{\neg U}(P)$ is a flat $\mathcal
O_W$-module. A coherent flat $\mathcal O_W$-module is necessarily a
locally free $\mathcal O_W$-module. Whence $\pi_*^{\neg U}(P)\in\cc
P(W,X)^{\neg U}$. Assertion (a) is proven.

Prove now assertion (b). Let $Z=X\setminus U$ be the complement of
$U$ and let $Z^{\prime}=X^{\prime}\setminus U^{\prime}$ and both
regarded as reduced schemes. The square $Q$ is an elementary
Nisnevich square, so $\pi$ induces a scheme isomorphism $Z^{\prime}
\to Z$. Let $y \in W \times X$. In this case
$\pi^{-1}(y)=\{y^{\prime}\}$ as sets for a unique point $y^{\prime}$
and $y^{\prime} \in W \times Z^{\prime}$. Moreover,
$\pi^{-1}(y)=Spec(k(y^{\prime}))$ as schemes and the map $\pi^{*}:
k(y) \to k(y^{\prime})$ is an isomporhism.

By (a) one has $\pi_*^{\neg U}(P)\in\cc P(W,X)^{\neg U}$. So the
composite morphism
$P\xrightarrow{adj}\pi^*\pi_*(P)\to\pi_*(\pi^*_{\neg U}(P))$ is a
morphism of coherent locally free $\mathcal O_W$-modules. Thus to
check that it is an isomorphism it suffices to check that for each
closed point $y \in Supp(P)$ the induced morphism $P(y) \to
\pi_*(\pi^*_{\neg U}(P))(y)$ is an isomorphism. Since $P\in\cc
P(W,X)^{\neg U}$ one has $y \in W \times Z$. In that case the
composite map
$$k(y) \xrightarrow{adj} \pi^*\pi_*(k(y))\to\pi_*(\pi^*_{\neg U}(k(y)))=\pi_*(k(y^{\prime}))=k(y^{\prime})$$
is an isomorphism. It follows that the induced morphism $P(y) \to
\pi_*(\pi^*_{\neg U}(P))(y)$ is an isomorphism as well. So the
morphism $P \to \pi_*(\pi^*_{\neg U}(P))$ from (b) is indeed an
isomorphism.

Finally prove (c). Since $Q \in P(W,X^{\prime})^{\neg U^{\prime}}$
then each closed point $y^{\prime} \in Supp(Q)$ belongs to $W \times
Z^{\prime}$, where $Z^{\prime} \subset X^{\prime}$ is from the proof
of (b). We already know that $\pi_*(Q) \in P(W,X)^{\neg U}$. Thus
(a) implies $\pi^*_{\neg U}(\pi_*(Q)) \in P(W,X^{\prime})^{\neg
U^{\prime}}$. So the map $\pi^*_{\neg
U}(\pi_*(Q))\to\pi^*\pi_*(Q)\xrightarrow{adj} Q$ is a morphism of
coherent locally free $\mathcal O_W$-modules. Therefore it suffices
to check that for each closed point $y^{\prime} \in Supp(Q)$ the map
\begin{equation}
\label{item_c} \pi^*_{\neg U}(\pi_*(Q))(y^{\prime})
\to\pi^*\pi_*(Q)(y^{\prime})\xrightarrow{adj} Q(y^{\prime})
\end{equation}
is an isomorphism. The latter follows from the fact that $y^{\prime}
\in W \times Z^{\prime}$ and the map
$$k(y^{\prime})= \pi^*_{\neg U}(\pi_*(k(y^{\prime}))) \to \pi^*(\pi_*(k(y^{\prime}))) \xrightarrow{adj} k(y^{\prime})$$
is identity. Assertion (c) is proven.
\end{proof}

The sublemma shows that
   $$\pi^*_{\neg U}:\cc P(W,X)^{\neg U}\leftrightarrows\cc P(W,X')^{\neg U'}:\pi_*^{\neg U}$$
are mutually inverse equivalences of categories. So we have shown
that $\cc O_K$ is Nisnevich excisive. The fact that $\cc
O_{K^\ps},\cc O_{K^{Gr}}$ are Nisnevich excisive is proved in a
similar way.

By above arguments it follows that for any elementary distinguished
square $Q$ the sequence of Nisnevich sheaves
   $$0\to\cc K_0(-,U')\to\cc K_0(-,U)\ps\cc K_0(-,X')\to\cc K_0(-,X)\to 0$$
is exact showing that $\cc O_{K_0}$ is Nisnevich excisive. The fact
that $\cc O_{K^\ps_0}$ is Nisnevich excisive is proved in a similar
way.

Now the sequence of Nisnevich sheaves
   $$0\to Cor(-,U')\to Cor(-,U)\ps Cor(-,X')\to Cor(-,X)\to 0$$
is exact by~\cite[4.3.9]{SV}. We conclude that $\cc O_{cor}$ is
Nisnevich excisive as well. The theorem is proved.
\end{proof}

\begin{cor}\label{qoqoqo}
$\cc O_{K^{\oplus}_0},\cc O_{K_0},\cc O_{cor}$ are symmetric
monoidal spectral categories, and hence motivically excisive.
\end{cor}

\begin{proof}
We have shown above that $\cc O_{cor}$ is symmetric monoidal (see
p.~\pageref{symmcor}). Let $X,Y,X',Y'$ be four smooth schemes and
let $P\in\cc P(X,Y),P'\in\cc P(X',Y')$. In this case, the external
tensor product $P\boxtimes P'$ is obviously finite and flat over
$X\times X'$. Thus, we get a bifunctor
   $$\boxtimes:\cc P(X,Y)\times\cc P(X',Y')\to\cc P(X\times X',Y\times Y'),$$
which is obviously additive and biexact. This gives a canonical
operation -- an external tensor product
   $$\boxtimes:K_0^\oplus(X,Y)\otimes K_0^\oplus(X',Y')\to K_0^\oplus(X\times X',Y\times Y')$$
and
   $$\boxtimes:K_0(X,Y)\otimes K_0(X',Y')\to K_0(X\times X',Y\times Y').$$
This external tensor product determines symmetric monoidal spectral
category structures for $\cc O_{K^{\oplus}_0}$ and $\cc O_{K_0}$.
Since $\cc O_{K^{\oplus}_0},\cc O_{K_0},\cc O_{cor}$ are Nisnevich
excisive by the preceding theorem, it follows from Lemma~\ref{qoqo}
that these are motivically excisive as well.
\end{proof}

\begin{rem}{\rm
We shall show below that $\cc O_{K^{Gr}}$, $\cc O_{K^\oplus}$, $\cc
O_K$ are motivically excisive spectral categories. For this we first
need to construct the bivariant motivic spectral sequence relating
bivariant $K$-theory to motivic cohomology.

}\end{rem}

The next theorem is the main result of the section.

\begin{thm}\label{modelmot}
Suppose $\cc O$ is a Nisnevich (respectively motivic) excisive
spectral category. Then the Nisnevich local (motivic) model
structure on $\Mod\cc O$ is cellular, proper, spectral and weakly
finitely generated. Moreover, a map of $\cc O$-modules is a weak
equivalence in the Nisnevich local (respectively motivic) model
structure \ifff it is a weak equivalence in the Nisnevich local
(respectively motivic) model structure on $Pre^\Sigma(Sm/F)$. If
$\cc O$ is a symmetric monoidal spectral category then each of the
model structures on $\Mod\cc O$ is symmetric monoidal with respect
to the smash product~\eqref{smash} of $\cc O$-modules.
\end{thm}

\begin{proof}
Since Bousfield localization respects cellularity and left
properness, then both model structures are cellular and left proper.
It follows from~\cite[3.2.13]{HSS} that every object of $\Mod\cc O$
is small. By Lemma~\ref{pqp} both model structures are weakly
finitely generated, because domains and codomains of morphisms from
$\cc N\cup\cc A$ are finitely presentable. These are plainly
spectral as well.

We are going to show the statement first for the Nisnevich local
model structure.

Let $J$ be a family of generating trivial stable projective
cofibrations for $Sp^\Sigma$. Notice that $J$ can be chosen in such
a way that domains and codomains of the maps in $J$ are finitely
presentable. Recall that the set of maps in $\Mod\cc O$
   $$\cc P_J=\{\cc O(-,X)\wedge sj\xrightarrow{\cc O(-,X)\wedge j}\cc \cc O(-,X)\wedge tj\mid j\in J,X\in Sm/F\}$$
is a family of generating trivial cofibrations for the stable
projective model structure on $\Mod\cc O$.

We set
   $$\wh{\cc N}_{\cc O}:=\{A\wedge\Delta[n]_+\coprod_{A\wedge\partial\Delta[n]_+}B\wedge\partial\Delta[n]_+
     \to B\wedge\Delta[n]_+\mid(A\to B)\in\cc N_{\cc O},n\geq 0\}.$$
Following terminology of~\cite[section~4.2]{Hir} an augmented family
of $\cc N_{\cc O}$-horns is the following family of trivial
cofibrations:
   $$\Lambda(\cc N_{\cc O})=\cc P_J\cup\wh{\cc N}_{\cc O}.$$
Observe that domains and codomains of the maps in $\Lambda(\cc
N_{\cc O})$ are finitely presentable. It can be proven similar
to~\cite[4.2]{Hov} that a map $f:A\to B$ is a fibration in the
Nisnevich local model structure with fibrant codomain if and only if
it has the right lifting property with respect to $\Lambda(\cc
N_{\cc O})$.

Consider the adjoint functors~\eqref{polezno}
   $$\xymatrix{{\Psi_*}:Pre^\Sigma(Sm/F)\ar@<0.5ex>[r]&\Mod\cc O:{\Psi^*}.\ar@<0.5ex>[l]}$$
We first observe that ${\Psi^*}$ takes every map in $\Lambda(\cc
N_{\cc O})$ to a weak equivalence in $Sp^{\Sigma}_{nis}(Pre(Sm/F))$.
Indeed, if $f\in\cc P_J$ then ${\Psi^*}(f)$ is a stable projective
weak equivalence, because ${\Psi^*}$ preserves stable projective
weak equivalences. Suppose $f\in\wh{\cc N}_{\cc O}$; then
$\Psi^*(f)$ is a weak equivalence in $Sp^{\Sigma}_{nis}(Pre(Sm/F))$,
because $\cc O$ is Nisnevich excisive by assumption.

We want next to check that $\Psi^*$ maps elements of $\Lambda(\cc
N_{\cc O})$-cell to weak equivalences in
$Sp^{\Sigma}_{nis}(Pre(Sm/F))$. Here $\Lambda(\cc N_{\cc O})$-cell
denotes the class of maps of sequential compositions of cobase
changes of coproducts of maps in $\Lambda(\cc N_{\cc O})$. Since
$\Psi^*$ preserves filtered colimits and weak equivalences are
closed under filtered colimits, it suffices to prove that $\Psi^*$
sends the cobase change of a map in $\Lambda(\cc N_{\cc O})$ to a
weak equivalence.

Clearly, $\Psi^*$ maps the cobase change of a map in $\cc P_J$ to a
stable weak equivalence in $Pre^\Sigma(Sm/F)$. Note that every
element in $\cc N_{\cc O}$ is a trivial cofibration in
$Sp^{\Sigma}_{nis}(Pre(Sm/F))$. Therefore every map in $\wh{\cc
N}_{\cc O}$ is a trivial cofibration in
$Sp^{\Sigma}_{nis}(Pre(Sm/F))$, and hence so is the cobase change of
every map in $\wh{\cc N}_{\cc O}$.

In order to show that $\Psi^*(f)$ is a weak equivalence in
$Sp^{\Sigma}_{nis}(Pre(Sm/F))$ for any weak equivalence $f$ in
$\Mod\cc O$, we use the small object argument
(see~\cite[10.5.16]{Hir} or~\cite[2.1.14]{H}). We construct a
fibrant replacement $\alpha:X\to L_{\cc N}X$ for $X\in\Mod\cc O$,
where $\alpha$ is the transfinite composition of a
$\aleph_0$-sequence
   $$X=E^0\lra{\alpha_0}E^1\lra{\alpha_1}E^2\lra{\alpha_2}\cdots$$
in which each $E^n\to E^{n+1}$ is constructed as follows. Let $S$ be
the set of all commutative squares
   $$\xymatrix{A\ar[r]\ar[d]_g&E^n\ar[d]\\
               B\ar[r]& {*} }$$
where $g\in\Lambda(\cc N_{\cc O})$. Then $\alpha_n$ a pushout
   $$\xymatrix{\coprod_{s\in S}A_s\ar[r]\ar[d]_{\coprod g_s}&E^n\ar[d]^{\alpha_n}\\
               \coprod_{s\in S}B_s\ar[r]&E^{n+1}}$$
This construction is functorial in $X$. We have verified that
$\alpha$ is a weak equivalence in $Sp^{\Sigma}_{nis}(Pre(Sm/F))$.

Now let $f:X\to Y$ be a weak equivalence in the Nisnevich local
model structure on $\Mod\cc O$. Then the diagram
   $$\xymatrix{\Psi^*(X)\ar[r]\ar[d]_{\Psi^*(f)}&\Psi^*(L_{\cc N}X)\ar[d]^{\Psi^*(L_{\cc N}f)}\\
               \Psi^*(Y)\ar[r]&\Psi^*(L_{\cc N}Y)}$$
is commutative. The horizontal arrows are weak equivalences in
$Sp^{\Sigma}_{nis}(Pre(Sm/F))$, the right arrow is a stable
projective weak equivalence. We infer that the left arrow is a weak
equivalence in $Sp^{\Sigma}_{nis}(Pre(Sm/F))$.

On the other hand, if $\Psi^*(f)$ is a weak equivalence in
$Sp^{\Sigma}_{nis}(Pre(Sm/F))$, then so is $\Psi^*(L_{\cc N}f)$. It
is, moreover, a stable projective weak equivalence
by~\cite[3.2.13]{Hir}, because $\Psi^*(L_{\cc N}X),\Psi^*(L_{\cc
N}Y)$ are fibrant in $Sp^{\Sigma}_{nis}(Pre(Sm/F))$. It follows that
$f$ is a weak equivalence in the Nisnevich local model structure on
$\Mod\cc O$. We have proved that a map of $\cc O$-modules is a weak
equivalence in the Nisnevich local model structure \ifff it is a
weak equivalence in $Sp^{\Sigma}_{nis}(Pre(Sm/F))$.

We claim that $\Psi^*$ respects fibrations. For this we shall apply
a theorem of Bousfield~\cite{Bou}. Consider the stable model
structure for $\cc O$-modules and a commutative diagram
   $$\xymatrix{L_{\cc N}V\ar[r]^{L_{\cc N}k}&L_{\cc N}X\\
               V\ar[u]\ar[d]\ar[r]^k&X\ar[u]\ar[r]^{\alpha_X}\ar[d]^f&L_{\cc N}X\ar[d]\\
               W\ar[d]\ar[r]^h&Y\ar[r]^{\alpha_Y}\ar[d]&L_{\cc N}Y\\
               L_{\cc N}W\ar[r]^{L_{\cc N}h}&L_{\cc N}Y}$$
with the central square pullback, $f$ a fibration between stably
fibrant objects and $\alpha_X,\alpha_Y,L_{\cc N}h$ stable projective
weak equivalences. Note that $\Psi^*(h)$ is a weak equivalence in
$Sp^{\Sigma}_{nis}(Pre(Sm/F))$. Since $\Psi^*:\Mod\cc O\to
Pre^\Sigma(Sm/F)$ is a right Quillen functor with respect to the
stable model structure, it follows that $\Psi^*(f)$ is a fibration
in $Pre^\Sigma(Sm/F)$. Since $\alpha_X,\alpha_Y$ are stable
projective weak equivalences and the stable model structure on
$Pre^\Sigma(Sm/F)$ is right proper, one easily sees that $\Psi^*$
takes the right square of the diagram to a homotopy pullback square.
Now~\cite[3.4.7]{Hir} implies $\Psi^*(f)$ is a fibration in
$Sp^{\Sigma}_{nis}(Pre(Sm/F))$. Since $Sp^{\Sigma}_{nis}(Pre(Sm/F))$
is right proper by Theorem~\ref{uhh}, we see that $\Psi^*(k)$ is a
weak equivalence in $Sp^{\Sigma}_{nis}(Pre(Sm/F))$. Thus
$\Psi^*(L_{\cc N}k)$ is a stable projective weak equivalence, and
hence so is $L_{\cc N}k$.

By~\cite[9.3, 9.7]{Bou} the following notions define a proper
simplicial model structure on $\Mod\cc O$: a morphism $f:X\to Y$ is
a cofibration if and only if it is a stable projective cofibration,
a weak equivalence if and only if $L_{\cc N}f:L_{\cc N}X\to L_{\cc
N}Y$ is a stable projective weak equivalence, and fibration if and
only if $f$ is a stable projective fibration and the commutative
square
   $$\xymatrix{X\ar[r]\ar[d]_f&L_{\cc N}X\ar[d]\\
               Y\ar[r]&L_{\cc N}Y}$$
is homotopy cartesian. This model structure plainly coincides with
the Nisnevich local model structure on $\Mod\cc O$, because
cofibrations and weak equivalences are the same. As a consequence,
$\Psi^*$ respects fibrations and is a right Quillen functor from $\Mod\cc O$ to
$Sp^{\Sigma}_{nis}(Pre(Sm/F))$. By~\cite[A.1.4]{Sch} if $\cc I$
($\cc J$) is a generating family of (trivial) cofibrations in
$Sp^{\Sigma}_{nis}(Pre(Sm/F))$, then $\Psi_*(\cc I)$ ($\Psi_*(\cc
J)$) is a generating family of (trivial) cofibrations.

Suppose $\cc O$ is a symmetric monoidal spectral category. Given two
cofibrations $i:A\to B,j:C\to D$, the map
   $$(i,j)_*:(B\wedge_{\cc O}C)\cup_{(A\wedge_{\cc O}C)}(A\wedge_{\cc O}D)\to B\wedge_{\cc O}D$$
is a cofibration, because cofibrations in the Nisnevich local and
stable projective model structures are the same and the latter model
structure is symmetric monoidal by Theorem~\ref{modelo}.

Let $i\in\Psi_*(\cc I)$ and $j\in\Psi_*(\cc J)$. Then $i=\Psi_*(i')$
and $j=\Psi_*(j')$ for some $(i':A'\to B')\in\cc I,(j':C'\to
D')\in\cc J$. The functor $\Psi_*$ is strong symmetric monoidal.
Therefore the map
   $$(i,j)_*:(B\wedge_{\cc O}C)\cup_{(A\wedge_{\cc O}C)}(A\wedge_{\cc O}D)\to B\wedge_{\cc O}D$$
is isomorphic to
   $$\Psi_*((i',j')_*):\Psi_*(B'\wedge_{\cc O_{naive}}C')\cup_{\Psi_*(A'\wedge_{\cc O_{naive}}C')}
     \Psi_*(A'\wedge_{\cc O_{naive}}D')\to\Psi_*(B'\wedge_{\cc O_{naive}}D'),$$
which is, in turn, isomorphic to
   $$\Psi_*((B'\wedge_{\cc O_{naive}}C')\cup_{(A'\wedge_{\cc O_{naive}}C')}
     (A'\wedge_{\cc O_{naive}}D'))\to\Psi_*(B'\wedge_{\cc O_{naive}}D').$$
The map
   \begin{equation*}\label{llaa}
    (B'\wedge_{\cc O_{naive}}C')\cup_{(A'\wedge_{\cc O_{naive}}C')}
    (A'\wedge_{\cc O_{naive}}D')\to(B'\wedge_{\cc O_{naive}}D')
   \end{equation*}
is a trivial cofibration in $Sp^\Sigma_{nis}(Pre(Sm/F))$ by
Theorem~\ref{uhh}. Since $\Psi_*$ respects trivial
cofibrations, then $(i,j)_*$ is a trivial cofibration in the
Nisnevich local model structure on $\Mod\cc O$. Therefore the
Nisnevich local model structure on $\Mod\cc O$ is symmetric monoidal
by~\cite[4.2.5]{H}.

To prove the statement for the motivic model structure on $\Mod\cc
O$, it is enough to verify that each map
   $$f:\cc O(-,X\times\bb A^1)\to\cc O(-,X),\quad X\in Sm/F,$$
is a weak equivalence in $Sp^{\Sigma}_{mot}(Pre(Sm/F))$, because the
rest of the proof is verified similar to the Nisnevich local model
structure. It is the case because we assume $\cc O$ to be
motivically excisive.
\end{proof}

\begin{cor}\label{vottak}
Suppose a spectral category $\cc O$ is Nisnevich excisive
(respectively motivically excisive). Then the pair of adjoint
fuctors
   $$\xymatrix{{\Psi_*}:Pre^\Sigma(Sm/F)\ar@<0.5ex>[r]&\Mod\cc O:{\Psi^*}\ar@<0.5ex>[l]}$$
induces a Quillen pair for the Nisnevich local projective
(respectively motivic) model structures on $Pre^\Sigma(Sm/F)$ and
$\Mod\cc O$. In particular, one has adjoint functors between
triangulated categories ${\Psi_*}:\shnis(F)\leftrightarrows\shnis\cc
O:{\Psi^*}$ (respectively
${\Psi_*}:\sheff(F)\leftrightarrows\sheff\cc O:{\Psi^*}$).
\end{cor}

To show that the map
   $$\cc O_K(-,X\times\bb A^1)\to\cc O_K(-,X)$$
is a motivic weak equivalence in $Pre^\Sigma(Sm/F)$ (and similarly
for $\cc O_{K^{\oplus}},\cc O_{K^{Gr}}$), one has to construct the
bivariant motivic spectral sequence. It will follow then from
Theorem~\ref{panin} that $\cc O_{K^{Gr}}$, $\cc O_{K^\oplus}$, $\cc
O_K$ are motivically excisive spectral categories.

\section{Bivariant motivic cohomology groups}

Consider a ringoid $\cc A$ and its Eilenberg-–Mac~Lane spectral
category $H\cc A$. We denote by $Ch(\cc A)$ the category of
unbounded chain complexes of $\cc A$-modules. By an $\cc A$-module
we just mean a contravariant additive functor from $\cc A$ to
abelian groups. For instance, if $\cc A$ is $Cor$ (respectively
$K_0^\oplus$ or $K_0$) then $\cc A$-modules are presheaves with
transfers (respectively with $K_0^\oplus$- or $K_0$-transfers).

By~\cite[section~B.1]{SS} there is a chain of Quillen equivalences
relating the category of $H\cc A$-modules $\Mod H\cc A$ with respect
to the stable projective model structure and $Ch(\cc A)$ with
respect to the usual model structure (the quasi-isomorphisms and
epimorphisms are weak equivalences and fibrations respectively):
   $$\xymatrix{\Mod H\cc A\ar@<1ex>[r]^U&\Nvmod H\cc A\ar@<1ex>[l]^L\ar@<1ex>[r]^(.6)\Lambda
               &Ch(\cc A)\ar@<1ex>[l]^(.4){\cc H}.}$$
Here $L,\Lambda$ are left adjoint and the intermediate model
category of naive $H\cc A$-modules is defined as follows.

\begin{defs}{\rm
Let $\cc O$ be a spectral category. A {\it naive\/} $\cc O$-module
$M$ consists of a collection $\{M(o)\}_{o\in\cc O}$ of $\bb Z_{\geq
0}$-graded, pointed simplicial sets together with associative and
unital action maps $M(o)_p\wedge O(o',o)_q\to M(o')_{p+q}$ for pairs
of objects $o,o'$ in $\cc O$ and for $p,q\geq 0$. A morphism of
naive $\cc O$-modules $M\to N$ consists of maps of graded spaces
$M(o)\to N(o)$ strictly compatible with the action of $\cc O$. We
denote the category of naive $\cc O$-modules by $\Nvmod\cc O$.

The {\it free\/} naive $\cc O$-module $F_o$ at an object $o\in\cc O$
is given by the graded spaces $F_o(o')=\cc O(o',o)$ with action maps
   $$F_o(o')_p\wedge\cc O(o'',o')_q=\cc O(o',o)_p\wedge\cc O(o'',o')_q\to\cc O(o'',o)_{p+q}=F_o(o'')_{p+q}$$
given by composition in $\cc O$.

}\end{defs}

One defines a model structure for naive $H\cc A$-modules as
follows~\cite[B.1.3]{SS}. A morphism of naive $H\cc A$-modules
$f:M\to N$ is a weak equivalence if it is an objectwise
$\pi_*$-isomorphism, i.e. if for all $a\in\cc A$ the map
$f(a):M(a)\to N(a)$ induces an isomorphism of stable homotopy
groups. The map $f$ is an objectwise stable fibration if each $f(a)$
is a stable fibration of spectra in the sense of~\cite[2.3]{BF}. A
morphism of naive $H\cc A$--modules is a cofibration if it has the
left lifting properties for maps which are objectwise
$\pi_*$-isomorphisms and objectwise stable firations.

The forgetful functor $U$ takes the free, genuine $H\cc A$-module to
the free, naive $H\cc A$-module. The left adjoint $L$ sends the
naive free modules $F_a$ to the genuine free modules. If we consider
the free $\cc A$-module $\cc A(-,a)$, as a complex in dimension 0,
then it is naturally isomorphic to $\Lambda(F_a)$
(see~\cite[section~B.1]{SS} for details).

The notions for fibrant naive $H\cc A$-modules and complexes of $\cc
A$-modules to be flasque are defined similar to $\cc O$-modules. We
say that a morphism $f:M\to N$ of $H\cc A$-modules (respectively
complexes of $\cc A$-modules) is a {\it Nisnevich local weak
equivalence\/} if for every flasque $H\cc A$-module $Q$ the morphism
$f^*:\Ho(\Nvmod H\cc A)(N[i],Q)\to\Ho(\Nvmod H\cc A)(M[i],Q)$
(respectively $f^*:\Ho(Ch\cc A)(N[i],Q)\to\Ho(Ch\cc A)(M[i],Q)$) is
an isomorphism for every integer $i$. Note that $f$ is a Nisnevich
local weak equivalence if and only if sheafification with respect to
Nisnevich topology of the graded morphism of graded presheaves
$\pi_*(f):\pi_*(M)\to\pi_*(N)$ (respectively $H^*(f):H^*(M)\to
H^*(N)$) is an isomorphism.

\begin{defs}{\rm
(1) A flasque $H\cc A$-module (respectively a complex of $\cc
A$-modules) $Q$ is said to be {\it $\bb A^1$-homotopy invariant\/}
if for every $X\in Sm/F$ the map $Q(X)\to Q(X\times\bb A^1)$ is a
$\pi_*$-isomorphism (respectively quasi-isomorphism).

(2) We say that a morphism $f:M\to N$ of $H\cc A$-modules
(respectively complexes of $\cc A$-modules) is a {\it motivic
equivalence\/} if for every $\bb A^1$-homotopy invariant $H\cc
A$-module $Q$ the morphism $f^*:\Ho(\Nvmod H\cc
A)(N[i],Q)\to\Ho(\Nvmod H\cc A)(M[i],Q)$ (respectively
$f^*:\Ho(Ch\cc A)(N[i],Q)\to\Ho(Ch\cc A)(M[i],Q)$) is an
isomorphism.

}\end{defs}

We define the Nisnevich local and motivic model structures for $H\cc
A$-modules or complexes of $\cc A$-modules as follows. The
cofibrations remain the same but the weak equivalences are the
Nisnevich local weak equivalences and motivic equivalences
respectively. Fibrations are defined by the corresponding lifting
property. It follows that the chain of adjoint functors
   $$\xymatrix{\Mod H\cc A\ar@<1ex>[r]^U&\Nvmod H\cc A\ar@<1ex>[l]^L\ar@<1ex>[r]^(.6)\Lambda
               &Ch(\cc A)\ar@<1ex>[l]^(.4){\cc H}}$$
yields Quillen equivalences between model categories with respect to
Nisnevich local and motivic model structures. Denote by
$D_{\nis}(\cc A)$ and $D_{\mot}(\cc A)$ the triangulated homotopy
categories of $Ch(\cc A)$ with respect to the Nisnevich local and
motivic model structures.

We get a pair of triangulated equivalences between triangulated
categories
   $$\xymatrix{\shnis H\cc A\ar@<1ex>[r]&D_{\nis}(\cc A).\ar@<1ex>[l]}$$

%

We want to give another description of $D_{\nis}(\cc A)$ for the
case when $\cc A$ is either $K^\oplus_0$ or $K_0$ or $Cor$. One
refers to $K^\oplus_0$- and $K_0$-modules as $K^\oplus_0$- and
$K_0$-presheaves respectively. $Cor$-modules are called in the
literature presheaves with transfers.

\begin{prop}[\cite{Sus,Voe1,Wlk}]\label{svw}
Let $\cc F$ be a $K_0^\oplus$-presheaf (respectively,
$K_0$-presheaf, presheaf with transfers). Then, the associated
Nisnevich sheaf $\cc F_{\nis}$ has a unique structure of a
$K_0^\oplus$-presheaf (respectively, $K_0$-presheaf, presheaf with
transfers) for which the canonical homomorphism $\cc F\to\cc
F_{\nis}$ is a homomorphism of $K_0^\oplus$-presheaves
(respectively, of $K_0$-presheaves, presheaves with transfers).
\end{prop}

Denote by $Sh({K_0^\oplus}),Sh(K_0),ShTr$ the categories of
Nisnevich $K_0^\oplus$-sheaves, $K_0$-sheaves and sheaves with
transfers respectively. Their derived categories are denoted by
$D(Sh({K_0^\oplus}))$, $D(Sh(K_0))$, $D(ShTr)$.

\begin{cor}\label{ivo}
$Sh({K_0^\oplus}),Sh(K_0),ShTr$ are Grothendieck categories, and
hence have enough injectives.
\end{cor}

\begin{proof}
We prove the claim for $Sh({K_0^\oplus})$, because the other two
cases are similarly checked. $Sh({K_0^\oplus})$ has filtered direct
limits which are exact, because this is the case for
$K_0^\oplus$-presheaves and for Nisnevich sheaves. So
$Sh({K_0^\oplus})$ satisfies axiom (Ab5). The category of
$K_0^\oplus$-presheaves is a Grothendieck category with
$\{K_0^\oplus(-,X)\}_{X\in Sm/F}$ the family of projective
generators. It follows that the family of sheaves
$\{K_0^\oplus(-,X)_{\nis}\}_{X\in Sm/F}$ is a family of generators
of $Sh({K_0^\oplus})$.
\end{proof}

Using the fact that Nisnevich local weak equivalences between
complexes of Nisnevich $K_0^\oplus$-sheaves (respectively
$K_0$-sheaves and sheaves with transfers) coincide with
quasi-isomorphisms between such complexes, we infer that the functor
of sheafification induces triangulated equivalences of triangulated
categories
   $$D_{\nis}({K_0^\oplus})\bl\sim\to D(Sh({K_0^\oplus})),\quad
     D_{\nis}(K_0)\bl\sim\to D(Sh(K_0)),\quad D_{\nis}(Cor)\bl\sim\to D(ShTr).$$
Note that $D_{\nis}({K_0^\oplus})$, $D_{\nis}(K_0)$ and
$D_{\nis}(Cor)$ are compactly generated triangulated categories with
compact generators given by representable presheaves, and hence so
are $D(Sh({K_0^\oplus})),D(Sh(K_0)),D(ShTr)$.

Given a smooth scheme $X$, we define $\bb Z_{K_0^\oplus}(X)$
(respectively $\bb Z_{K_0}(X)$ and $\bb Z_{tr}(X)$) as the complex
having $K_0^\oplus(-,X)_{\nis}$ (respectively $K_0(-,X)_{\nis}$ and
$Cor(-,X)$) in degree zero and zero in other degrees. Here
$K_0^\oplus(-,X)_{\nis}$ stands for the Nisnevich sheaf associated
to the presheaf $U\mapsto K_0^\oplus(U,X)$.

\begin{defs}{\rm
(1) We say that a complex of $K_0^\oplus$-sheaves (respectively
$K_0$-sheaves and sheaves with transfers) $Q$ is {\it $\bb
A^1$-local\/} if for every scheme $X\in Sm/F$ and every integer $n$
the natural map
   $$D(Sh({K_0^\oplus}))(\bb Z_{K_0^\oplus}(X)[n],Q)\to D(Sh({K_0^\oplus}))(\bb Z_{K_0^\oplus}(X\times\bb A^1)[n],Q)$$
(respectively the maps $D(Sh({K_0}))(\bb Z_{K_0}(X)[n],Q)\to
D(Sh({K_0}))(\bb Z_{K_0}(X\times\bb A^1)[n],Q)$ and $D(ShTr)(\bb
Z_{tr}(X)[n],Q)\to D(ShTr)(\bb Z_{tr}(X\times\bb A^1)[n],Q)$) is an
isomorphism.

(2) A morphism $M\to N$ of complexes of $K_0^\oplus$-sheaves
(respectively $K_0$-sheaves and sheaves with transfers) is called an
{\it $\bb A^1$-weak equivalence\/} if for every $\bb A^1$-local
complex $Q$ the map
   $$D(Sh({K_0^\oplus}))(N,Q)\to D(Sh({K_0^\oplus}))(M,Q)$$
(respectively the maps $D(Sh({K_0}))(N,Q)\to D(Sh({K_0}))(M,Q)$ and
$D(ShTr)(N,Q)\to D(ShTr)(\bb Z_{tr}(M,Q)$) is an isomorphism.

(3) The {\it $\bb A^1$-derived category\/} of $K_0^\oplus$-sheaves
(respectively $K_0$-sheaves and sheaves with transfers) is the one,
obtained from $D(Sh({K_0^\oplus}))$ (respectively $D(Sh({K_0}))$ and
$D(ShTr)$) by inverting the $\bb A^1$-weak equivalences. The
corresponding $\bb A^1$-derived categories will be denoted by
$D_{\bb A^1}(Sh({K_0^\oplus}))$, $D_{\bb A^1}(Sh({K_0}))$ and
$D_{\bb A^1}(ShTr)$ respectively.

}\end{defs}

By the general localization theory of compactly generated
triangulated categories (see, e.g.,~\cite{N}) $D_{\bb
A^1}(Sh({K_0^\oplus}))$ (respectively $D_{\bb A^1}(Sh({K_0}))$ and
$D_{\bb A^1}(ShTr)$) is the localization of $D(Sh({K_0^\oplus}))$
(respectively $D(Sh({K_0}))$ and $D(ShTr)$) with respect to the
localizing subcategory generated by cochain complexes of the form
   $$\bb Z_{K_0^\oplus}(X\times\bb A^1)\to\bb Z_{K_0^\oplus}(X)$$
(respectively $\bb Z_{K_0}(X\times\bb A^1)\to\bb Z_{K_0}(X)$ and
$\bb Z_{tr}(X\times\bb A^1)\to\bb Z_{tr}(X)$).

$D_{\bb A^1}(Sh({K_0^\oplus}))$ (respectively $D_{\bb
A^1}(Sh({K_0}))$ and $D_{\bb A^1}(ShTr)$) can be identified with the
full subcategory of $D(Sh({K_0^\oplus}))$ of $\bb A^1$-local
complexes. The inclusion functor
   $$D_{\bb A^1}(Sh({K_0^\oplus}))\to D(Sh({K_0^\oplus}))$$
admits a left adjoint
   $$L_{\bb A^1}:D(Sh({K_0^\oplus}))\to D_{\bb A^1}(Sh({K_0^\oplus}))$$
which is also called the {\it $\bb A^1$-localization functor}. The
same $\bb A^1$-localization functor exists for $D_{\bb
A^1}(Sh({K_0}))$ and $D_{\bb A^1}(ShTr)$.

The above arguments may be summarized as follows.

\begin{prop}\label{docum}
There are natural equivalences
   $$D_{\mot}({K_0^\oplus})\bl\sim\to D_{\bb A^1}(Sh({K_0^\oplus})),\quad
     D_{\mot}(K_0)\bl\sim\to D_{\bb A^1}(Sh(K_0)),\quad D_{\mot}(Cor)\bl\sim\to D_{\bb A^1}(ShTr).$$
of triangulated categories.
\end{prop}

Recall that a sheaf $\cc F$ of abelian groups in the Nisnevich
topology on $Sm/F$ is {\it strictly $\bb A^1$-invariant\/} if for
any $X\in Sm/F$, the canonical morphism
   $$H^*_{\nis}(X,\cc F)\to H^*_{\nis}(X\times\bb A^1,\cc F)$$
is an isomorphism. It follows from~\cite[6.2.7]{Mor} that a complex
of $K_0^\oplus$-sheaves (respectively $K_0$-sheaves and sheaves with
transfers) $Q$ is $\bb A^1$-local if and only if each cohomology
sheaf $H^n(Q)$ is strictly $\bb A^1$-invariant.

Recall that for any presheaf of abelian groups $\cc F$ on $Sm/F$ we
get a simplicial presheaf $C_n(\cc F)$, by setting $C_n(\cc
F)(U)=\cc F(U\times\Delta^n)$. We shall write $C^*(\cc F)$ to denote
the corresponding cochain complex (of degree $+1$) of abelian
presheaves. Namely,
   $$C^i(\cc F):=C_{-i}(\cc F).$$
By~\cite{Sus,Voe1,Wlk} strictly $\bb A^1$-invariant $K_0^\oplus$-,
$K_0$-sheaves and sheaves with transfers coincide with $\bb
A^1$-invariant ones whenever the field $F$ is perfect. Therefore
over perfect fields the categories $D_{\bb A^1}(Sh({K_0^\oplus}))$,
$D_{\bb A^1}(Sh({K_0}))$ and $D_{\bb A^1}(ShTr)$ can be identified
with the full subcategories of complexes having $\bb A^1$-invariant
cohomology sheaves. Moreover, the functor $L_{\bb A^1}$ equals the
functor $C^*$ in all three cases. A detailed proof for sheaves with
transfers is given in~\cite{Voe1}. The triangulated category of
Voevodsky $DM_-^{eff}$~\cite{Voe1} is a full subcategory of $D_{\bb
A^1}(ShTr)$ in this case.

We want to introduce bivariant motivic cohomology for smooth
schemes, but first recall some facts for motivic cohomology.  Let
$X\in Sm/F$, $\star\in\{K^\oplus_0,K_0,tr\}$, and $n\geq 0$; we
define the Nisnevich sheaf $\bb Z_\star(X)(\bb G_m^{\wedge n})$ as
follows. Let $\cc D_n$ be the sum of images of homomorphisms
   $$\bb Z_\star(X\times\bb G_m^{\times n-1})\to\bb Z_\star(X\times\bb G_m^{\times n})$$
induced by the embeddings of the form
   $$(a_1,\ldots,a_{n-1})\in\bb G_m^{\times n-1}\mapsto(a_1,\ldots,1,\ldots,a_{n-1})\in\bb G_m^{\times n}.$$
The sheaf $\bb Z_\star(X)(\bb G_m^{\wedge n})$ is, by definition,
$\bb Z_\star(X)(\bb G_m^{\times n})/\cc D_n$. In what follows we
shall denote the sheaf $\bb Z_\star(pt)(\bb G_m^{\wedge n})$ by $\bb
Z_\star(\bb G_m^{\wedge n})$.

\begin{defs}{\rm
The {\it $K_0^\oplus$-motive $M_{K_0^\oplus}(X)$\/} of a smooth
scheme $X$ over $F$ (respectively the {\it $K_0$-motive\/}
$M_{K_0}(X)$ and the {\it motive $M_{tr}(X)$}) is the image of $\bb
Z_{K_0^\oplus}(X)$ in $D_{\bb A^1}(Sh({K_0^\oplus}))$ (respectively
the image of $\bb Z_{K_0}(X)$ in $D_{\bb A^1}(Sh(K_0))$ and the
image of $\bb Z_{tr}(X)$ in $D_{\bb A^1}(ShTr)$). If
$\star\in\{K_0^\oplus,K_0,tr\}$ then by $M_\star(X)(n)$, $n\geq 0$,
we denote the corresponding image of $\bb Z_\star(X)(\bb G_m^{\wedge
n})[-n]$.

}\end{defs}

Let $X\in Sm/F$; the {\it complex $\bb Z_\star(X)(n)$ of weight
$n$\/} on $Sm/F$ is the complex $C^*\bb Z_\star(X)(\bb G_m^{\wedge
n})[-n]$, where the degree shift refers to the cohomological
indexing of complexes. If $X=pt$ we shall just write $\bb
Z_\star(n)$ for the motivic complex dropping $X$ from notation.

\begin{defs}{\rm
For smooth schemes $U,X\in Sm/F$ we define their {\it bivariant
motivic cohomology groups $H^{i,n}(U,X,\bb Z)$\/} as
$H^i_{\nis}(U,\bb Z(X)_{K_0^{\oplus}}(n))$. We shall write
$H^{i,n}_{\cc M}(U,\bb Z)$ to denote $H^i_{\nis}(U,\bb
Z_{K_0^{\oplus}}(n))$.

}\end{defs}

By a theorem of Suslin~\cite{Sus} for any $n\geq 0$ and any field
$F$, the canonical homomorphism of complexes of Nisnevich sheaves
   $$f_n:\bb Z_{K_0^{\oplus}}(n)\to\bb Z_{tr}(n)$$
is a quasi-isomorphism. Hence, for any smooth scheme $X\in Sm/F$,
cohomology $H^*_{\nis}(X,\bb Z_{K_0^{\oplus}}(n))$ coincides with
motivic cohomology $H^*_{\nis}(X,\bb Z_{tr}(n))$ of
Suslin--Voevod\-sky~\cite{SV1}. By~\cite{Voe2} for any field $F$,
any smooth scheme $X$ over $F$ and any $i,n\in\bb Z$, there is a
natural isomorphism
   $$H^i_{\nis}(X,\bb Z_{tr}(n))\cong CH^n(X,2n-i),$$
where the right hand side groups are higher Chow groups. Since both
groups are homotopy invariant, then $H^i_{\nis}(X,\bb
Z_{K_0^{\oplus}}(n))$ are homotopy invariant
$K_0^{\oplus}$-presheaves.

\begin{prop}\label{interny}
For any Nisnevich $K_0^\oplus$-sheaf $\cc F$ we have natural
identifications
   $$\Ext^i_{Sh(K_0^\oplus)}(\bb Z_{K_0^\oplus}(X),\cc F)=H^i_{\nis}(X,\cc F).$$
In particular $\Ext^i_{Sh(K_0^\oplus)}(\bb Z_{K_0^\oplus}(X),-)=0$
for $i>\dim X$.
\end{prop}

\begin{proof}
Since the category $Sh(K_0^\oplus)$ has sufficiently many injective
objects by Corollary~\ref{ivo} and for any $K_0^\oplus$-sheaf $\cc
G$ one has $\Hom(K_0^\oplus(-,X)_{\nis},\cc G)=\cc G(X)$ we only
have to show that for any injective $K_0^\oplus$-sheaf $\cc I$ one
has $H^i_{\nis}(X,\cc I)=0$ for $i>0$.

It follows from the proof of Theorem~\ref{panin} that for any
elementary distinguished square~\eqref{squareQ} the sequence of
$K_0^\oplus$-sheaves
   $$0\to K_0^\oplus(-,U')_{\nis}\to K_0^\oplus(-,U)_{\nis}\oplus K_0^\oplus(-,X')_{\nis}
     \to K_0^\oplus(-,X)_{\nis}\to 0$$
is exact. Therefore the sequence of abelian groups
   $$0\to\cc I(X)\to\cc I(U)\oplus\cc I(X')\to\cc I(U')\to 0$$
is exact, because $\cc I$ is injective. Let $H\cc I$ be the
Eilenberg--Mac~Lane sheaf of $S^1$-spectra associated with $\cc I$
(see~\cite[p.~23]{Mor}). It follows that $H\cc I$ is Nisnevich
excisive.

By~\cite[3.2.3]{Mor} there is a natural isomorphism
   $$H^n_{\nis}(X,\cc I)\to[(X_+),H\cc I[n]],\quad n\in\bb Z,$$
where the right hand side is the Hom-set in the stable homotopy
category of $S^1$-spectra on $Sm/F$ (see~\cite{Mor}). It follows
from~\cite[3.1.7]{Mor} that $[(X_+),H\cc I[n]]$ is isomorphic to
$\pi_{-n}(H\cc I(X))$ which is zero for $n\not=0$.
\end{proof}

\begin{cor}\label{pampam}
Let $X$ be a smooth scheme over a field $F$ and $C$ be a complex of
Nisnevich $K_0^\oplus$-sheaves bounded above. Then for any $i\in\bb
Z$ there is a canonical isomorphism
   $$\Hom_{D(Sh(K_0^\oplus))}(\bb Z_{K_0^\oplus}(X),C[i])\cong H^i_{\nis}(X,C).$$
\end{cor}

\begin{proof}
The proof is like that of~\cite[1.8]{SV1} (one should use the
preceding proposition as well).
\end{proof}

\begin{cor}\label{tututu}
For every $n\geq 0$ the complex $\bb Z_{K_0^\oplus}(n)$ is $\bb
A^1$-local in $D(Sh(K_0^\oplus))$.
\end{cor}

\begin{proof}
This follows from the previous corollary and the fact that
$H^i_{\nis}(X,\bb Z_{K_0^\oplus}(n))$ are homotopy invariant
$K_0^\oplus$-presheaves.
\end{proof}

\begin{thm}\label{agree}
Let $F$ be any field, then for every $X\in Sm/F$ there is a natural
isomorphism
   $$H^{i,n}_{\cc M}(X,\bb Z)\cong D_{\bb A^1}(Sh(K_0^\oplus))(M_{K_0^\oplus}(X),M_{K_0^\oplus}(pt)(n)[i]).$$
If the field $F$ is perfect then there is also a natural isomorphism
   $$H^{i,n}(U,X,\bb Z)\cong D_{\bb A^1}(Sh(K_0^\oplus))(M_{K_0^\oplus}(U),M_{K_0^\oplus}(X)(n)[i])$$
for any $U,X\in Sm/F$.
\end{thm}

\begin{proof}
The preceding corollary implies $\bb Z_{K_0^\oplus}(n)$ is $\bb
A^1$-local in $D(Sh(K_0^\oplus))$. It is proved similar
to~\cite[section~3.2]{Voe1} that for any Nisnevich
$K_0^\oplus$-sheaf $\cc F$, the natural morphism of complexes
   $$\cc F\to C^*(\cc F)$$
is an $\bb A^1$-weak equivalence. By Corollary~\ref{pampam} there is
a canonical isomorphism
   $$\Hom_{D(Sh(K_0^\oplus))}(\bb Z_{K_0^\oplus}(X),\bb Z_{K_0^\oplus}(n)[i])\cong H^{i,n}_{\cc M}(X,\bb Z).$$
for any $i\in\bb Z$. On the other hand one has,
   $$\Hom_{D(Sh(K_0^\oplus))}(\bb Z_{K_0^\oplus}(X),\bb Z_{K_0^\oplus}(n)[i])\cong
     \Hom_{D(Sh(K_0^\oplus))}(\bb Z_{K_0^\oplus}(X)(0),\bb Z_{K_0^\oplus}(n)[i]).$$
Thus,
   $$H^{i,n}_{\cc M}(X,\bb Z)\cong \Hom_{D(Sh(K_0^\oplus))}(\bb Z_{K_0^\oplus}(X)(0),\bb
   Z_{K_0^\oplus}(n)[i])\cong D_{\bb A^1}(Sh(K_0^\oplus))(M_{K_0^\oplus}(X),M_{K_0^\oplus}(pt)(n)[i]).$$
is an isomorphism for all $i\in\bb Z$.

Assume now that $F$ is perfect. Then all cohomology sheaves for $\bb
Z_{K_0^\oplus}(X)(n)$ are homotopy invariant, hence strictly
homotopy invariant because $F$ is perfect (see above). We see that
each complex $\bb Z_{K_0^\oplus}(X)(n)$ is $\bb A^1$-local. The
second assertion is now checked similar to the first one.
\end{proof}

\section{The Grayson tower}

Recall that $Ord$ denotes the category of finite nonempty ordered
sets, and for each $d\geq 0$ we introduce the object
$[d]=\{0<1<\cdots<d\}$ of $Ord$. Given $A,B\in Ord$ we let $AB\in
Ord$ denote the ordered set obtained by concatenating $A$ and $B$,
with the elements of $A$ smaller than the elements of $B$. Given a
simplicial set $Y$ with base point $y_0\in Y_0$, the natural
inclusion maps $A\to AB\leftarrow B$ provide natural face maps
$Y(A)\to Y(AB)\leftarrow Y(B)$. Let $PY$ be the simplicial path
space of edges in $Y$ with initial endpoint at $y_0$; it can be
defined for $A\in Ord$ by
   $$(PY)(A)=\lo(\{y_0\}\to Y([0])\leftarrow Y([0]A)).$$
The space $|PY|$ is contractible. The face maps $Y([0]A)\to Y(A)$
provide a projection map $PY\to Y$. We define
   $$\omega Y=\lo(PY\to Y\leftarrow PY).$$
The commutative square
   $$\xymatrix{\omega Y\ar[r]\ar[d]&PY\ar[d]\\
               PY\ar[r]& Y}$$
together with the contractibility of the space $|PY|$ provides a
natural map $|wY|\to\Omega|Y|$.

For instance, let $Y$ be the nerve of a category $\cc B$ with
$b_0\in\Ob\cc B$ the base point. A vertex of $\omega Y$ is nothing
more than a pair of morphisms in $\cc B$
   $$b_0\to b_1\leftarrow b_0,\quad b_1\in\Ob\cc B.$$
An edge of $\omega Y$ is a pair of commutative triangles
   $$\xymatrix{b_0\ar[r]^u\ar[dr]&b_1\ar[d]&b_0\ar[l]_v\ar[dl]\\
               &b_2}$$
The map $|\omega Y|\to\Omega|Y|$ is uniquely determined by a map
$|\omega Y|\wedge S^1\to|Y|$. The latter map comes from two maps
$|\omega Y|\times|\Delta^1|\to|Y|$, which are uniquely determined by
two functors of categories $U,V:\cc B\times\{0\to 1\}\to\cc B$.
These functors can be thought of as two natural transformations from
the constant functor sending all objects of $\cc B$ to $b_0$ to the
identity functor on $\cc B$. The natural transformations are defined
by the arrows $u,v$ on objects and the two triangles define them on
morphisms. A vertex of $\omega Y$ yields a loop in $|Y|$ which
starts at $b_0$, follows the edge $u$ to $b_1$, and returns along
the other edge $v$ to $b_0$, and an edge of $\omega Y$ yields a
homotopy between two such loops.

If $\cc M$ is an exact category, and $S$ is the $S$-construction of
Waldhausen~\cite{Wal} then $\omega S\cc M$ is the simplicial set
$G\cc M$ of Gillet--Grayson~\cite{GiGr} and $|\omega S\cc
M|\to\Omega|SM|$ is a homotopy equivalence as it was shown
in~\cite{GiGr}. Also, if $\cc M$ is an additive category, then
$|\omega S^\oplus\cc M|\to\Omega|S^\oplus\cc M|$ is a homotopy
equivalence~\cite{Gr}, where $S^\oplus$ is the Grayson
$S^\oplus$-construction. Note that the latter equivalence is
functorial in $\cc M$.

Given an additive category $\cc M$, Quillen defines a new category
$S^{-1}S\cc M$ whose objects are pairs $(A,B)$ of objects of $\cc
M$. A morphism $(A,B)\to(C,D)$ in $S^{-1}S\cc M$ is given by a pair
of split monomorphisms
   $$f:A\bl{\twoheadleftarrow}\rightarrowtail C,\quad g:B\bl\twoheadleftarrow\rightarrowtail D$$
together with an isomorphism $h:\coker f\to\coker g$. By a split
monomorphism we mean a monomorphism together with a chosen
splitting. The nerve of the category $S^{-1}S\cc M$ which is also
denoted by $S^{-1}S\cc M$ is homotopy equivalent to Quillen's
$K$-theory space of $\cc M$ by~\cite{Gr1}. There is an obvious map
   $$u:\omega S^\oplus\cc M\to S^{-1}S\cc M.$$
It is a homotopy equivalence by~\cite[section~4]{Gr}.

We can regard $\omega S^\oplus\cc M$ and $S^{-1}S\cc M$ as
simplicial additive categories. The equivalence $u$ above can be
extended to a stable equivalence of symmetric spectra
   $$u:K^{Gr}(\omega S^\oplus\cc M)\to K^{Gr}(S^{-1}S\cc M).$$

Though the space $|\omega S^{-1}S\cc M|$ is not the loop space of
$S^{-1}S\cc M$~\cite[section~4]{Gr}, it is worth to consider
$|\omega S^{-1}S\cc M|$ by replacing $\cc M$ by a simplicial
additive category over a connected simplicial ring.

Suppose now that $X$ is a pointed simplicial space. We let $I_*X_d$
denote the connected component of $X_d$ containing the base point.
Grayson~\cite{Gr} has shown that there is an essential obstruction
to the natural map $|d\mapsto\Omega X_d|\to\Omega|X|$ being an
equivalence. The following theorem says what this obstruction is
when $X$ is a simplicial group-like $H$-space.

\begin{thm}[Grayson~\cite{Gr}]\label{grayson}
If $X$ is a simplicial group-like $H$-space then the natural
sequence
\begin{equation}\label{eqgr}
  |d\mapsto\Omega X_d|\to\Omega|X|\to\Omega|d\mapsto\pi_0(X_d)|
\end{equation}
is a fibration sequence. Moreover, the map
   $$|X|\to|d\mapsto\pi_0(X_d)|$$
induces an isomorphism on $\pi_0$, and its homotopy fiber is
connected.
\end{thm}

The standard diagonalization technique allows us to
generalize~\eqref{eqgr} to multisimplicial spaces. For example, if
$X$ is a bisimplicial group-like $H$-space, then
   $$|(d,e)\mapsto\Omega X_{d,e}|\to\Omega|X|\to\Omega|(d,e)\mapsto\pi_0(X_{d,e})|$$
is a fibration sequence.

If $R$ is a simplicial ring, and $\cc M$ is a simplicial additive
category, we say that $\cc M$ is {\it $R$-linear\/} if, for each
$d\geq 0$, $M_d$ is an $R_d$-linear category, and for each map
$\phi:[e]\to[d]$, each $r\in R_d$, and each arrow $f$ in $\cc M_d$,
we have the equation $\phi^*(rf)=\phi^*(r)\phi^*(f)$. A typical
example of a contractible ring is $F[\Delta]$ whose $n$-simplices
are defined as
   $$F[\Delta]_n=F[x_0,\ldots,x_n]/(x_0+\cdots+x_n-1).$$
By $\Delta^{\cdot}$ we denote the cosimplicial affine scheme
$\spec(F[\Delta])$.

\begin{thm}[Grayson~\cite{Gr}]\label{grayson2}
If $R$ is a contractible simplicial ring and $\cc M$ is a simplicial
$R$-linear additive category, then the natural map
   $$|d\mapsto|\omega S^{-1}S\cc M_d||\to|d\mapsto\Omega|S^{-1}S\cc M_d||$$
is a homotopy equivalence of spaces.
\end{thm}

Here is some convenient notation for cubes. We let $[1]$ denote the
ordered set $\{0<1\}$ regarded as a category, and we use $\epsilon$
as notation for an object of $[1]$. By an $n$-dimensional cube in a
category $\cc C$ we will mean a functor from $[1]^n$ to $\cc C$. An
object $C$ in $\cc C$ gives a 0-dimensional cube denoted by $[C]$,
and an arrow $C\to C'$ in $\cc C$ gives a 1-dimensional cube denoted
by $[C\to C']$. If the category $\cc C$ has products, we may define
an external product of cubes as follows. Given an $n$-dimensional
cube $X$ and an $n'$-dimensional cube $Y$ in $\cc C$, we let
$X\boxtimes Y$ denote the $n+n'$-dimensional cube defined by
$(X\boxtimes
Y)(\epsilon_1,\ldots,\epsilon_{n+n'})=X(\epsilon_l,\ldots,\epsilon_n)\times
Y(\epsilon_{n+1},\ldots,\epsilon_{n+n'})$. Let $\bb G_m^{\wedge n}$
denote the external product of $n$ copies of $[1\to \bb G_m]$. For
example, $\bb G^{\wedge 2}_m$ is the square of schemes
   $$\xymatrix{\spec F\ar[r]\ar[d]&\bb G_m\ar[d]\\
               \bb G_m\ar[r]&\bb G_m\times\bb G_m.}$$

Denote by $\cc M\langle\bb G_m^{\times n}\rangle$ the exact category
where an object is a tuple $(P,\theta_1,\ldots,\theta_n)$ consisting
of an object $P$ of $\cc M$ and commuting automorphisms
$\theta_1,\ldots,\theta_n\in Aut(P)$. Note that $\cc M\langle\bb
G_m^{\times n}\rangle=(\cc M\langle\bb
G_m^{\times(n-1)}\rangle)\langle\bb G_m\rangle$. For instance, if
$\cc M=\cc P(U,X)$, $U,X\in Sm/F$, then we can identify $\cc
P(U,X)\langle\bb G_m^{\times n}\rangle$ with $\cc P(U,X\times\bb
G_m^{\times n})$. The cube of affine schemes $\bb G_m^{\wedge n}$
gives rise to a cube of exact categories $\cc M\langle\bb
G_m^{\wedge n}\rangle$ with vertices being $\cc M\langle\bb
G_m^{\times k}\rangle$, $0\leq k\leq n$. The edges of the cube are
given by the natural exact functors $i_s:\cc M\langle\bb
G_m^{\times(k-1)}\rangle\to\cc M\langle\bb G_m^{\times k}\rangle$
defined as
   $$(P,(\theta_1,\ldots,\theta_{k-1}))\longmapsto(P,(\theta_1,\ldots,1,\ldots,\theta_{k-1})),$$
where 1 is the $s$th coordinate.

In~\cite[\S4]{Gr2} is presented a construction called $C$ which can
be applied to a cube of exact categories to convert it into a
multisimplicial exact category, the $K$-theory of which serves as
the iterated cofiber space of the corresponding cube of $K$-theory
spaces/spectra.

Given an exact category $\cc M$ with a chosen zero object 0 and an
ordered set $A$, we call a functor $F:\Ar(A)\to\cc M$ exact if
$F(i,i)=0$ for all $i$, and $0\to F(i,j)\to F(i,k)\to F(j,k)\to 0$
is {\it exact\/} for all $i\leq j\leq k$. The set of such exact
functors is denoted by $Exact(\Ar(A),\cc M)$.

Now let $L$ be a symbol, and consider $\{L\}$ to be an ordered set.
Given an $n$-dimensional cube of exact categories $\cc M$, we define
an $n$-fold multisimplicial exact category $C\cc M$ as a functor
from $(Ord^n)^{\op}$ to the category of exact categories by letting
$C\cc M(A_1,\ldots,A_n)$ be the set
   $$Exact([\Ar(A_1)\to\Ar(\{L\}A_1)]\boxtimes\cdots\boxtimes[\Ar(A_n)\to\Ar(\{L\}A_n)],\cc M)$$
of multi-exact natural transformations. When $n=0$, we may identify
$C\cc M$ with $\cc M$. In the case $n=1$, $C\cc M$ is the same as a
construction of Waldhausen~\cite{Wal} denoted $S.(\cc M_0\to\cc
M_1)$. We define $S.\cc M$ to be $S.C\cc M$, the result of applying
the $S.$-construction of Waldhausen degreewise. The construction
$S.\cc M$ is a $n+1$-fold multisimplicial set. $K$-theory of $C\cc
M$ serves as the iterated cofiber space/spectrum of the
corresponding cube of $K$-theory spaces/spectra.

Given an $n$-dimensional cube of additive categories $\cc M$, we
define an $n$-fold multisimplicial additive category $C^\oplus\cc M$
as a functor from $(Ord^n)^{\op}$ to the category of exact
categories by letting $C^\oplus\cc M(A_1,\ldots,A_n)$ be the set
   $$Add([Sub(A_1)\to Sub(\{L\}A_1)]\boxtimes\cdots\boxtimes[Sub(A_n)\to Sub(\{L\}A_n)],\cc M)$$
of multi-additive natural transformations. When $n=0$, we may
identify $C^\oplus\cc M$ with $\cc M$. We define $S^\oplus\cc M$ to
be $S^\oplus C^\oplus\cc M$, the result of applying the
$S^\oplus$-construction of Grayson degreewise. It is an $n+1$-fold
multisimplicial set (see~\cite{Gr} for details). Grayson's
$K$-theory $K^{Gr}(C^\oplus\cc M)$ of $C^\oplus\cc M$ (respectively
Waldhausen's $K$-theory $K(C\cc M)$ of $C\cc M$) serves as the
iterated cofiber space/spectrum of the corresponding cube of
Grayson's (Waldhausen's) $K$-theory spaces/spectra. It is easy to
see that
   $$K^{Gr}_0(C^\oplus\cc M\langle\bb G_m^{\wedge n}\rangle)=K^{Gr}_0(\cc M\langle\bb G_m^{\times n}\rangle)/
     \sum^n_{k=1}(i_k)_*(K^{Gr}_0(\cc M\langle\bb G_m^{\times (n-1)}\rangle)).$$

If $\cc M$ is an additive category, Grayson showed~\cite{Gr} that
the category $\omega S^{-1}S\cc M$ is equivalent to the category
$\cc C$ where an object is any pair $(P,\beta)$ with $P\in\cc M$ and
$\beta\in Aut(P)$, and an arrow $(P,\beta)\to(Q,\gamma)$ is any
split monomorphism $P\bl{\twoheadleftarrow}\rightarrowtail Q$ with
respect to which one has the equation $\gamma=\beta\oplus 1$. The
objects of $\cc C$ are themselves the objects of the exact category
$\cc M\langle\bb G_m\rangle$, the arrows of $\cc C$ are also the
objects of an exact category, and indeed, the nerve of $\cc C$ can
be interpreted as a simplicial exact category.

There is a map $C^\oplus\cc M\langle\bb G_m^{\wedge 1}\rangle\to\cc
C$ which amounts to forgetting some choices of cokernels, so the map
is a homotopy equivalence. Thus we have a homotopy equivalence of
symmetric $K$-theory spectra
   $$v:K^{Gr}(C^\oplus\cc M\langle\bb G_m^{\wedge 1}\rangle)\lra{\sim}K^{Gr}(\omega S^{-1}S\cc M).$$

Consider the category of topological symmetric spectra
$TopSp^\Sigma$ (see~\cite[section~I.1]{Sch}). We can apply adjoint
functors ``geometric realization", denoted by $|-|$, and ``singular
complex", denoted by $\cc S$, levelwise to go back and forth between
simplicial and topological symmetric spectra
   \begin{equation}\label{dwdw}
    |-|:Sp^\Sigma\rightleftarrows TopSp^\Sigma:\cc S.
   \end{equation}

\begin{rem}{\rm
By the standard abuse of notation $|-|$ denotes both the functor
from $Sp^\Sigma$ to $TopSp^\Sigma$ and the realization functor from
simplicial spectra to spectra. It will always be clear from the
context which of either meanings is used.

}\end{rem}

We have a zig-zag of maps in $TopSp^\Sigma$ between semistable
symmetric spectra
   \begin{gather*}
       |K^{Gr}(C^\oplus\cc M\langle\bb G_m^{\wedge 1}\rangle)|\xrightarrow{|v|}|K^{Gr}(\omega S^{-1}S\cc M)|\xrightarrow{w}\Omega|K^{Gr}(S^{-1}S\cc M)|
       \xleftarrow{\Omega|u|}\\ \Omega|K^{Gr}(\omega S^\oplus\cc M)|\xrightarrow{\sim}\Omega(\Omega|K^{Gr}(\cc M)|[1])|\xleftarrow{\sim}\Omega|K^{Gr}(\cc M)|.
   \end{gather*}
All maps except $w$ are weak equivalences of ordinary spectra. The
zig-zag yields an arrow in $\Ho(Sp^\Sigma)$
   \begin{equation}\label{prpr}
    \gamma:K^{Gr}(C^\oplus\cc M\langle\bb G_m^{\wedge 1}\rangle)\to\Omega K^{Gr}(\cc M).
   \end{equation}
Observe that the arrow is functorial in $\cc M$. If $\cc M$ is a
simplicial $R$-linear additive category over a contractible
simplicial ring, then
   $$|d\mapsto K^{Gr}(C^\oplus\cc M_d\langle\bb G_m^{\wedge 1}\rangle)|\lra{\gamma}|d\mapsto\Omega K^{Gr}(\cc M_d)|$$
is an isomorphism in $\Ho(Sp^\Sigma)$ by Theorem~\ref{grayson2}. In
fact, $\gamma$ is an isomorphism in the homotopy category $\Ho(Sp)$
of ordinary spectra $Sp$, because all arrows in the zig-zag
producing $\gamma$ are weak equivalences in $Sp$.

Thus Theorem~\ref{grayson} implies that the following sequence is a
triangle in $\Ho(Sp^\Sigma)$
   $$|d\mapsto K^{Gr}(C^\oplus\cc M_d\langle\bb G_m^{\wedge 1}\rangle)|\lra{}\Omega|d\mapsto K^{Gr}(\cc M_d)|
     \lra{}\Omega|d\mapsto EM(K_0(\cc M_d))|\lra{+},$$
where $EM(K_0(\cc M_d))$ is the Eilenberg--Mac~Lane spectrum of $\cc
M_d$ (see Appendix). This triangle yields a triangle
   $$S^1\wedge|d\mapsto K^{Gr}(C^\oplus\cc M_d\langle\bb G_m^{\wedge 1}\rangle)|\lra{}|d\mapsto K^{Gr}(\cc M_d)|
     \lra{}|d\mapsto EM(K_0(\cc M_d))|\lra{+}.$$
We call it the {\it Grayson triangle}. It produces more generally
triangles
   $$S^1\wedge|d\mapsto K^{Gr}(C^\oplus\cc M_d\langle\bb G_m^{\wedge n+1}\rangle)|\to|d\mapsto K^{Gr}(\cc M_d)\langle\bb G_m^{\wedge n}\rangle|
     \to|d\mapsto EM(K_0(\cc M_d\langle\bb G_m^{\wedge n}\rangle))|\bl+\to.$$

In what follows we denote by $\cc K_0^\oplus$ and $\cc K_0$ the
Eilenberg--Mac~Lane spectral categories associated with ringoids
$K_0^\oplus$ and $K_0$ on $Sm/F$. We refer the reader to Appendix to
read about Eilenberg--Mac~Lane spectral categories. There are
canonical morphisms of spectral categories
   $$\cc O_{K^\oplus}\to\cc K_0^\oplus\leftarrow\cc O_{K_0^\oplus}$$
and
   $$\cc O_K\to\cc K_0\leftarrow\cc O_{K_0},$$
where the arrows on the right are equivalences of spectral
categories. It follows from~\cite[A.1.1]{SS} that restriction and
extension of scalars functors induce spectral Quillen equivalences
of modules
   $$\Mod\cc K_0^\oplus\rightleftarrows\Mod\cc O_{K_0^\oplus},\quad\Mod\cc K_0\rightleftarrows\Mod\cc O_{K_0}.$$

Consider the case when each $\cc M_d=\cc P'(U\times\Delta^d,X)$,
$U,X\in Sm/F$. Let $S_{K^{Gr}}(X)(n)$ (respectively $S_{\cc
K_0^\oplus}(X)(n)$) denote the presheaf $S^n\wedge|d,U\mapsto
K^{Gr}(C^\oplus\cc P'(U\times\Delta^d,X)\langle\bb G_m^{\wedge
n}\rangle)|$ (respectively the presheaf $S^n\wedge|d,U\mapsto
EM(K_0^\oplus(C^\oplus\cc P'(U\times\Delta^d,X))\langle\bb
G_m^{\wedge n}\rangle)|$). Note that $S_{\cc K_0^\oplus}(X)(n)$ is
an $\cc O_{K^{Gr}}$-module by means of the natural morphism of
spectral categories
   $$\cc O_{K^{Gr}}\to\cc K_0^\oplus.$$
We have a canonical morphism of presheaves
   $$g_q:S_{K^{Gr}}(X)(q)\to S_{\cc K_0^\oplus}(X)(q),\quad q\geq 0.$$
The Grayson triangle produces triangles in the homotopy category
$\Ho(\Mod\cc O_{naive})$ of $\cc O_{naive}$-modules regarded as a
model category with respect to the stable projective model structure
   $$S_{K^{Gr}}(X)(q+1)\xrightarrow{f_{q+1}}S_{K^{Gr}}(X)(q)\lra{g_q}S_{\cc K_0^\oplus}(X)(q)\bl+\to.$$

Recall that $Sp^{\bb N}_{nis}(Pre(Sm/F))$ stands for the model
category of $S^1$-spectra on pointed simplicial presheaves with
respect to the Nisnevich local model structure. Below we shall need
the following couple of lemmas.

\begin{lem}[\cite{FS}]\label{kino}
Consider a sequence of maps of spectra
   $$\cdots\xrightarrow{f_{q+2}}X_{q+1}\xrightarrow{f_{q+1}}X_q\lra{f_q}\cdots\lra{f_1}X_0=X$$
Assume further that for each $q$ we are given a map of spectra
$p_q:X_q\to B_q$ such that the composition $X_{q+1}\to X_q\to B_q$
is trivial and the associated map from $X_{q+1}$ to the homotopy
fiber of $X_q\to B_q$ is a weak equivalence. Assume further that for
each $i\geq 0$ there exists $n\geq 0$ such that $X_q$ is
$i$-connected for $q\geq n$. In this case there exists a strongly
convergent spectral sequence
   $$E^2_{pq}=\pi_{p+q}(B_q)\Longrightarrow \pi_{p+q}(X).$$
\end{lem}

We shall write $[E,L]$ to denote $\Ho(Sp^{\bb
N}_{nis}(Pre(Sm/F)))(E,L)$ for any two pre\-sheaves of spectra
$E,L$.

\begin{lem}\label{motyga}
Let $E$ be a presheaf of spectra and $U\in Sm/F$ of Krull dimension
$d$. Assume further that for each $i<q$ the $i$-th stable homotopy
Nisnevich sheaf $\pi_i(E)$ of $E$ is zero. Then for all $n<q-d$ one
has:
   $$[U_+[n],E]=0.$$
\end{lem}

\begin{proof}
This follows from~\cite[3.3.3]{Mor}.
\end{proof}

\begin{defs}{\rm
(1) Given a smooth scheme $X$ over $F$, the {\it Grayson tower\/} is
the sequence of maps in $\Ho(\Mod\cc O_{naive})$:
   $$\cdots\xrightarrow{f_{q+2}}S_{K^{Gr}}(X)(q+1)\xrightarrow{f_{q+1}}S_{K^{Gr}}(X)(q)\lra{f_q}\cdots\lra{f_1}S_{K^{Gr}}(X)(0).$$
By construction, the Grayson tower naturally produces a tower in
$\Ho(Sp^{\bb N}_{nis}(Pre(Sm/F)))$.

(2) For any $U,X\in Sm/F$ the {\it bivariant $K$-theory groups\/}
are defined as
   $$K_i(U,X)=[U_+[i],S_{K^{Gr}}(X)(0)],\quad i\in\bb Z.$$

}\end{defs}

\begin{lem}\label{tpy}
For any $U\in Sm/F$ and any integer $i$ there is a natural
isomorphism $K_i(U)\cong K_i(U,pt)$.
\end{lem}

\begin{proof}
This follows from Thomason's theorem~\cite{T} stating that algebraic
$K$-theory satisfies Nisnevich descent and the fact that $K(U)$ is
homotopy invariant.
\end{proof}

Denote by $Sp^{\bb N}_{nis,J}(Pre(Sm/F))$ (``$J$" for Jardine) the
stable model category of presheaves of ordinary spectra
corresponding to the Nisnevich local injective model structure on
$Pre(Sm/F)$ (see~\cite{Jar3}). Note that the identity functor
   $$Sp^{\bb N}_{nis}(Pre(Sm/F))\to Sp^{\bb N}_{nis,J}(Pre(Sm/F))$$
induces a left Quillen equivalence.

\begin{prop}\label{pumpum}
For any $U,X\in Sm/F$ and any $p,q$ there is a natural isomorphism
   $$H^{p,q}(U,X,\bb Z)\cong[U_+,S_{\cc K_0^\oplus}(X)(q)[p-2q]].$$
\end{prop}

\begin{proof}
Since the map of spectral categories $\cc O_{K_0^\oplus}\to\cc
K_0^\oplus$ is a levelwise weak equivalence, it is enough to show
the assertion for the $\cc O_{K_0^\oplus}$-module
   $$S'_{K_0^\oplus}(X)(q)=S^q\wedge|d,U\mapsto H(K_0^\oplus(C^\oplus\cc P'(U\times\Delta^d,X)\langle\bb G_m^{\wedge q}\rangle))|.$$

Denote by $Sp^{\bb N}_{nis}(s\Mod_{\bb Z})$ the model category of
simplicial $\bb Z$-module spectra in the sense of
Jardine~\cite{Jar4}. Every such spectrum consists of a sequence of
simplicial presheaves of abelian groups $A^n$, $n\geq 0$, together
with simplicial homomorphisms $A^n\otimes S^1\to A^{n+1}$, which are
also called bonding maps. The suspension spectrum
$\Sigma_{S^1}^\infty\cc X$ in $Sp^{\bb N}_{nis}(s\Mod_{\bb Z})$ of a
simplicial presheaf of abelian groups $\cc X$ is defined in the
usual way. We set,
   $$S''_{K_0^\oplus}(X)(q):=S^q\otimes|d,U\mapsto\Sigma_{S^1}^\infty(K_0^\oplus(C^\oplus\cc P'(U\times\Delta^d,X)\langle\bb G_m^{\wedge q}\rangle))|.$$

There is a natural forgetful functor
   $$U:Sp^{\bb N}_{nis}(s\Mod_{\bb Z})\to Sp^{\bb N}_{nis,J}(Pre(Sm/F)).$$
A map $f:A\to B$ of simplicial $\bb Z$-module spectra is a weak
equivalence if the underlying map of presheaves of spectra $Uf:UA\to
UB$ is a weak equivalence in $Sp^{\bb N}_{nis}(Pre(Sm/F))$. In fact,
$U$ is a right Quillen functor whose left adjoint is denoted by $V$.
The spectra $U(S''_{K_0^\oplus}(X)(q))$ and $S'_{K_0^\oplus}(X)(q)$
are stably equivalent by~\cite[5.5]{Jar4}. Therefore there is an
isomorphism
   $$\Ho(Sp^{\bb N}_{nis,J}(Pre(Sm/F)))(U_+,S'_{K_0^\oplus}(X)(q))\cong\Ho(Sp^{\bb N}_{nis}(s\Mod_{\bb Z}))(U_+,S''_{K_0^\oplus}(X)(q)).$$

It follows from~\cite{Jar4} that there are pairs of functors of
triangulated categories
   $$\xymatrix{\Ho(Sp^{\bb N}_{nis,J}(Pre(Sm/F)))\ar@<1ex>[r]^(.55)V&\Ho(Sp^{\bb N}_{nis}(s\Mod_{\bb Z}))\ar@<1ex>[l]^(.45)U\ar@<1ex>[r]^(.55)N
               &D(Sh(Sm/F))\ar@<1ex>[l]^(.45){\Gamma}.}$$
Here $D(Sh(Sm/F))$ is the derived category of Nisnevich sheaves,
$N,\Gamma$ are mutually inverse equivalences, and $V,N$ are left
adjoint. Moreover, $N$ takes the cofibrant presheaf of spectra
$S''_{K_0^\oplus}(X)(q)\in Sp^{\bb N}_{nis}(s\Mod_{\bb Z})$ to $\bb
Z_{K_0^\oplus}(X)(q)[2q]$. Our assertion now follows from the fact
that $Sp^{\bb N}_{nis}(Pre(Sm/F)),Sp^{\bb N}_{nis,J}(Pre(Sm/F))$ are
Quillen equivalent by means of the identity functor.
\end{proof}

We shall say that a presheaf of spectra $E$ is {\it $n$-connected\/}
if for all $k\leq n$ its $k$-th homotopy sheaf $\pi_k(E)$ vanishes.
We are now in a position to prove the main result of this section.

\begin{thm}\label{spektralka}
For every $q\geq 0$ and every $X\in Sm/F$ the sequence of maps
   $$S_{K^{Gr}}(X)(q+1)\xrightarrow{f_{q+1}}S_{K^{Gr}}(X)(q)\lra{g_{q}}S_{\cc K_0^\oplus}(X)(q)\lra{+}$$
is a triangle in $\Ho(Sp^{\bb N}_{nis}(Pre(Sm/F)))$ and the Grayson
tower produces a strongly convergent spectral sequence which we
shall call the {\em bivariant motivic spectral sequence},
   $$E_2^{pq}=H^{p-q,-q}(U,X,\bb Z)\Longrightarrow K_{-p-q}(U,X),\quad U\in Sm/F.$$
Moreover, if $X=pt$ and $U$ is the spectrum of a smooth Henselian
$F$-algebra then Grayson's spectral sequence~\cite{Gr} takes the
form
   $$E_2^{pq}=H_{\cc M}^{p-q,-q}(U,\bb Z)\Longrightarrow K_{-p-q}(U).$$
\end{thm}

\begin{proof}
The presheaf of spectra $S(X)(q)$, $q\geq 0$, is $(q-1)$-connected.
Lemmas~\ref{kino}-\ref{motyga} and Proposition~\ref{pumpum} imply
that the Grayson tower produces a strongly convergent spectral
sequence
   $$E_2^{pq}=H^{p-q,-q}(U,X,\bb Z)\Longrightarrow K_{-p-q}(U,X).$$

If $U$ is the spectrum of a smooth Henselian $F$-algebra then for
every presheaf of spectra $E$ one has:
   $$[U_+[n],E]\cong\pi_n(E(U)),\quad n\in\bb Z.$$
For $X=pt$, evaluation of the Grayson tower at $U$ is isomorphic in
$\Ho(Sp)$ to the tower constructed by Walker~\cite{Wlk}. The latter
tower produces a spectral sequence which agrees with Grayson's
spectral sequence~\cite{Gr}.
\end{proof}

We want to construct an $\bb A^1$-local counterpart for the
bivariant motivic spectral sequence. We denote by
   $$L_{\bb A^1}:Sp^{\bb N}_{nis}(Pre(Sm/F))\to Sp^{\bb N}_{nis}(Pre(Sm/F))$$
the $\bb A^1$-localization functor of Morel~\cite{Mor}. For any
presheaf $E$ of spectra and any integer $n$, the sheaves
   $$\pi_n^{\bb A^1}(E):=\pi_n(L_{\bb A^1}(E))$$
are strictly $\bb A^1$-invariant~\cite{Mor}.

\begin{defs}{\rm
(1) Given $X\in Sm/F$ and $q\geq 0$, the {\it Grayson $\cc
O_{K^{Gr}}$-module $G(X)(q)$ of weight $q$\/} is the $\cc
O_{K^{Gr}}$-module $S^q\wedge K^{Gr}(C^\oplus\cc P'(-,X)\langle\bb
G_m^{\wedge q}\rangle)$. We shall also write $G_0(X)(q)$ to denote
the $\cc K_0^\oplus$-module $S^q\wedge EM(K_0^\oplus(C^\oplus\cc
P'(-,X)\langle\bb G_m^{\wedge q}\rangle))$.

(2) For any $U,X\in Sm/F$ the {\it $\bb A^1$-local bivariant
$K$-theory groups\/} are defined as
   $$K_i^{\bb A^1}(U,X):=[U_+[i],L_{\bb A^1}(G(X)(0))],\quad i\in\bb Z.$$

(3) The {\it $\bb A^1$-local bivariant motivic cohomology groups\/}
are defined as
   $$H_{\bb A^1}^{p,q}(U,X,\bb Z):=[U_+,L_{\bb A^1}(G_0(X)(q))[p-2q]].$$

}\end{defs}

\begin{lem}\label{motyga2}
Let $E$ be a connective presheaf of spectra and $U\in Sm/F$ of Krull
dimension $d$. Then the group
   $$[U_+,L_{\bb A^1}(E)[n]]$$
vanishes for $n>d$.
\end{lem}

\begin{proof}
This follows from~\cite[4.3.1]{Mor}.
\end{proof}

The map~\eqref{prpr} yields a map in the stable homotopy category
$\Ho(\Mod\cc O_{K^{Gr}})$ of $\cc O_{K^{Gr}}$-modules
   $$f_q:G(X)(q)\to G(X)(q-1),\quad q>0.$$
In turn, one has a natural map of $\cc O_{K^{Gr}}$-modules
   $$g_q:G(X)(q)\to G_0(X)(q),\quad q\geq 0.$$

\begin{defs}{\rm
(1) Given a smooth scheme $X$ over $F$, the {\it Grayson tower of
Grayson's modules\/} is the sequence of maps in $\Ho(\Mod\cc
O_{K^{Gr}})$:
   $$\cdots\xrightarrow{f_{q+2}}G(X)(q+1)\xrightarrow{f_{q+1}}G(X)(q)\lra{f_q}\cdots\lra{f_1}G(X)(0).$$

(2) We say that a presheaf of spectra $E$ is {\it $\bb A^1$-local\/}
if for every scheme $U\in Sm/F$ and every integer $n$ the natural
map
   $$[(U_+)[n],E]\to[((U\times\bb A^1)_+)[n],E]$$
is an isomorphism.

}\end{defs}

We denote by $Sp^{\bb N}_{mot}(Pre(Sm/F))$ the model category of
$S^1$-spectra associated to the projective motivic model structure
on $Pre(Sm/F)$.

\begin{thm}\label{spektralka2}
For every $q\geq 0$ and every $X\in Sm/F$ the sequence of maps
   $$G(X)(q+1)\xrightarrow{f_{q+1}}G(X)(q)\lra{g_{q}}G_0(X)(q)\lra{+}$$
is a triangle in $\Ho(Sp^{\bb N}_{mot}(Pre(Sm/F)))$ and the Grayson
tower of Grayson's modules produces a strongly convergent spectral
sequence which we shall call the {\em $\bb A^1$-local bivariant
motivic spectral sequence},
   $$E_2^{pq}=H_{\bb A^1}^{p-q,-q}(U,X,\bb Z)\Longrightarrow K^{\bb A^1}_{-p-q}(U,X),\quad U\in Sm/F.$$
Assume further that one of the following conditions is satisfied:
\begin{enumerate}
\item $X=pt$;

\item the field $F$ is perfect.
\end{enumerate}
Then $H_{\bb A^1}^{p,q}(U,X,\bb Z)$ agree with bivariant motivic
cohomology groups $H^{p,q}(U,X,\bb Z)$ and the $\bb A^1$-local
bivariant motivic spectral sequence coincides with the bivariant
motivic spectral sequence of Theorem~\ref{spektralka}.
\end{thm}

\begin{proof}
We have a commutative diagram in the stable homotopy category of
presheaves of spectra
   $$\xymatrix{G(X)(q+1)\ar[r]^{f_{q+1}}\ar[d]&G(X)(q)\ar[r]^{g_{q}}\ar[d]&G_0(X)(q)\ar[d]\ar[r]^(.7){+}&{}\\
               S_{K^{Gr}}(X)(q+1)\ar[r]^(.53){f_{q+1}}&S_{K^{Gr}}(X)(q)\ar[r]^{g_q}&S_{\cc K_0^\oplus}(X)(q)\ar[r]^(.7){+}&{}}$$
with vertical arrows level $\bb A^1$-equivalences. Since the lower
sequence is a triangle in $\Ho(Sp^{\bb N}_{mot}(Pre(Sm/F)))$ then so
is the upper one. The presheaf of spectra $G(X)(q)$, $q\geq 0$, is
$(q-1)$-connected. Lemmas~\ref{kino} and~\ref{motyga2} imply that
the Grayson tower of Grayson's modules produces a strongly
convergent spectral sequence
   $$E_2^{pq}=H_{\bb A^1}^{p-q,-q}(U,X,\bb Z)\Longrightarrow K^{\bb A^1}_{-p-q}(U,X).$$

Assume now that $X=pt$. By Corollaries~\ref{pampam}-\ref{tututu} and
Proposition~\ref{pumpum} each $S_{\cc K_0^\oplus}(X)(q)$, $q\geq 0$,
is $\bb A^1$-local. In turn, if $F$ is perfect but $X$ is any smooth
scheme then Theorem~\ref{agree} and Proposition~\ref{pumpum} imply
each $S_{\cc K_0^\oplus}(X)(q)$, $q\geq 0$, is $\bb A^1$-local. In
both cases therefore $H_{\bb A^1}^{p,q}(U,X,\bb Z)$ agree with
bivariant motivic cohomology groups $H^{p,q}(U,X,\bb Z)$. The fact
that the $\bb A^1$-local bivariant motivic spectral sequence
coincides with the bivariant motivic spectral sequence of
Theorem~\ref{spektralka} is now obvious.
\end{proof}

\begin{cor}\label{trew}
There is a natural isomorphism of abelian groups
   $$K^{\bb A^1}_i(U,pt)\cong K_i(U)$$
for any $i\in\bb Z$. If $F$ is a perfect field, then there is also a
natural isomorphism
   $$K^{\bb A^1}_i(U,X)\cong K_i(U,X)$$
for any $i\in\bb Z$ and $U,X\in Sm/F$.
\end{cor}

\begin{proof}
This follows from Lemma~\ref{tpy} and the preceding theorem.
\end{proof}

We conclude the section by noting that the presheaves of $K$-groups
   $$K_i^{Gr}(\cc P'(-,Y))=\pi_i(\cc O_{K^{Gr}}(-,Y))$$
are different from both $K_i(-,Y)$ and $K_i^{\bb A^1}(-,Y)$ in
general. Indeed, suppose $F$ is perfect. Then the preceding
corollary implies the presheaves $K_i(-,Y)$ and $K_i^{\bb A^1}(-,Y)$
are isomorphic. Let $SH\cc O_{K^{Gr}}$ be the homotopy category of
$\cc O_{K^{Gr}}$-modules with respect to the stable projective model
structure (see Theorem~\ref{modelo}). It is a compactly generated
triangulated category. If the presheaves $K_i^{Gr}(\cc P'(-,Y))$
were isomorphic to $K_i^{\bb A^1}(-,Y)$ then we would have that the
natural map
   $${SH\cc O_{K^{Gr}}}(\cc O_{K^{Gr}}(-,X),\cc O_{K^{Gr}}(-,Y)[n])\to
     {SH\cc O_{K^{Gr}}}(\cc O_{K^{Gr}}(-,X\times\bb A^1),\cc O_{K^{Gr}}(-,Y)[n])$$
is an isomorphism for every $X,Y\in Sm/k$ and $n\in\bb Z$. Since
$\{\cc O_{K^{Gr}}(-,Y)[n]\}_{n\in\bb Z,Y\in Sm/k}$ is a family of
compact generators for $SH\cc O_{K^{Gr}}$, it would follow that the
natural map
   $$\cc O_{K^{Gr}}(-,X\times\bb A^1)\to\cc O_{K^{Gr}}(-,X)$$
is an isomorphism in $SH\cc O_{K^{Gr}}$ what is not the case.

\section{$K$-motives}

We use the $\bb A^1$-local bivariant motivic spectral sequence to
prove the following

\begin{thm}\label{nakonets}
The spectral categories $\cc O_{K^{Gr}},\cc O_{K^\oplus},\cc O_K$
are motivically excisive.
\end{thm}

\begin{proof}
The spectral categories $\cc O_{K^{Gr}},\cc O_{K^\oplus},\cc O_K$
are Nisnevich excisive by Theorem~\ref{panin}. We first check that
the natural map
   $$\cc O_{K^{Gr}}(-,X\times\bb A^1)\to\cc O_{K^{Gr}}(-,X),\quad X\in Sm/F,$$
is a motivic weak equivalence in $Pre^\Sigma(Sm/F)$.
By~\cite[4.34]{Jar2} it is enough to show that the map is a motivic
weak equivalence of presheaves of ordinary spectra. Since $\cc
O_{K_0^\oplus}$ is motivically excisive, then so is $\cc K_0^\oplus$
because these are equivalent spectral categories. In fact, the
natural map
   $$\cc K_0^\oplus(-,X\times\bb A^1)\to\cc K_0^\oplus(-,X),\quad X\in Sm/F,$$
is a level motivic weak equivalence of presheaves of ordinary
spectra.

Observe that the exact category $\cc P(U,X)\langle\bb G_m^{\times
q}\rangle$, $U,X\in Sm/F$, can be identified with $\cc
P(U,X\times\bb G_m^{\times q})$. It follows that the map
   $$\cc K_0^\oplus(-,X\times\bb A^1\times\bb G_m^{\times q})\to\cc K_0^\oplus(-,X\times\bb G_m^{\times q})$$
is a motivic weak equivalence of presheaves of ordinary spectra, and
hence so is the map
   $$G_0(X\times\bb A^1)(q)\to G_0(X)(q).$$
It induces an isomorphism of $\bb A^1$-local bivariant motivic
cohomology groups
   $$H_{\bb A^1}^{p,q}(U,X\times\bb A^1,\bb Z)\lra\cong H_{\bb A^1}^{p,q}(U,X,\bb Z).$$
We infer that the natural map of Grayson towers
   $$\xymatrix{\cdots\ar[r]^(.3){f_{q+1}}&G(X\times\bb A^1)(q)\ar[r]^(.45){f_q}\ar[d]&G(X\times\bb A^1)(q-1)\ar[d]\ar[r]^(.73){f_{q-1}}&\cdots\\
               \cdots\ar[r]^(.4){f_{q+1}}&G(X)(q)\ar[r]^{f_q}&G(X)(q-1)\ar[r]^(.7){f_{q-1}}&\cdots}$$
produces an isomorphism of $\bb A^1$-local bivariant motivic
spectral sequences, and hence each map
   $$K^{\bb A^1}_{i}(U,X\times\bb A^1)\to K^{\bb A^1}_{i}(U,X)$$
is an isomorphism. We see that
   $$L_{\bb A^1}(G(X\times\bb A^1)(0))\to L_{\bb A^1}(G(X)(0))$$
is a motivic weak equivalence, and hence so is
      $$\cc O_{K^{Gr}}(-,X\times\bb A^1)=G(X\times\bb A^1)(0)\to\cc O_{K^{Gr}}(-,X)=G(X)(0).$$
So $\cc O_{K^{Gr}}$ is motivically excisive. Since $\cc O_{K^{Gr}}$
and $\cc O_{K^{\oplus}}$ are equivalent spectral categories, then
also $\cc O_{K^\oplus}$ is motivically excisive.

It follows from~\cite[10.5]{Gr} that the natural map of spectra
   $$|n\mapsto K^\oplus(\cc P'(X\times\Delta^n,Y))|\to|n\mapsto K(\cc P'(X\times\Delta^n,Y))|$$
is a stable equivalence of spectra. Therefore,
   $$|n\mapsto\cc O_{K^{\oplus}}(X\times\Delta^n,Y)|
     \to|n\mapsto\cc O_K(X\times\Delta^n,Y)|$$
is a stable equivalence of spectra. For any presheaf $\cc F:Sm/F\to
SSets$, the natural map $\cc F\to Sing(\cc F)$ is a motivic
equivalence. Here $Sing:Pre(Sm/F)\to Pre(Sm/F)$ is the singular
functor (see, e.g.,~\cite[p.~542]{Jar2}). We conclude that for any
$X\in Sm/F$ the map of presheaves of spectra
      $$\cc O_{K^{\oplus}}(-,X)\to\cc O_K(-,X)$$
is a motivic equivalence, and therefore it is a motivic equivalence
of presheaves of symmetric spectra by~\cite[4.34]{Jar2}. Since $\cc
O_{K^{\oplus}}$ is motivically excisive, then so is $\cc O_K$.
\end{proof}

We say that a presheaf of symmetric spectra $E\in Pre^\Sigma(Sm/F)$
is {\it semistable\/} if $E(U)$ is a semistable symmetric spectrum
for every $U\in Sm/F$. We remark that all presheaves of symmetric
spectra we work with in practice like $\cc O_{K^{Gr}}(-,X)$, $\cc
O_{K^\oplus}(-,X)$, $\cc O_K(-,X)$, $\cc O_{K_0^\oplus}(-,X)$, $\cc
O_{K_0}(-,X)$, $\cc K_0^\oplus(-,X)$, $\cc K_0(-,X)$ are semistable.

\begin{lem}\label{qwer}
Let $E$ be a semistable presheaf of symmetric spectra and $U\in
Sm/F$. Then there are natural isomorphisms
   $$[U_+[n],E]\cong\shnis(F)(U_+[n],E),\quad[U_+[n],L_{\bb A^1}(E)]\cong\sheff(F)(U_+[n],E)$$
for all integers $n$.
\end{lem}

\begin{proof}
Straightforward.
\end{proof}

\begin{defs}\label{refer}{\rm
(1) The {\it $K^\oplus$-motive $M_{K^\oplus}(X)$\/} of a smooth
scheme $X$ over $F$ (respectively the {\it $K$-, $K_{Gr}$-, $\cc
K_0^\oplus$-, $\cc K_0$-motives\/} $M_K(X)$, $M_{K_{Gr}}(X)$,
$M_{\cc K_0^\oplus}(X)$, $M_{\cc K_0}(X)$) is the image of $\cc
O_{K^\oplus}(-,X)$ in $\sheff\cc O_{K^\oplus}$ (respectively the
corresponding images of $\cc O_K(-,X)$, $\cc O_{K_{Gr}}(-,X)$, $\cc
O_{\cc K_0^\oplus}(-,X)$, $\cc O_{\cc K_0}(-,X)$ in $\sheff\cc O_K$,
$\sheff\cc O_{K_{Gr}}$, $\sheff\cc O_{\cc K_0^\oplus}$, $\sheff\cc
O_{\cc K_0}$).

(2) The image of each Grayson module $G(X)(q)$ (respectively
$G_0(X)(q)$), $q\geq 0$, in $\sheff\cc O_{K^{Gr}}$ (respectively in
$\sheff\cc K_0^{\oplus}$) will be denoted by $M_{K^{Gr}}(X)(q)$
(respectively $M_{\cc K_0^{\oplus}}(X)(q)$).

}\end{defs}

\begin{lem}\label{aqe}
For any $X,Y\in Sm/F$ there are canonical isomorphisms:
\begin{gather*} K_i^{\bb A^1}(X,Y)\cong\sheff\cc
O_{K^{Gr}}(M_{K^{Gr}}(X)[i],M_{K^{Gr}}(Y))\cong\\ \sheff\cc
O_{K^\oplus}(M_{K^\oplus}(X)[i],M_{K^\oplus}(Y))\cong\sheff\cc
O_{K}(M_{K}(X)[i],M_{K}(Y))\end{gather*}
and
\begin{gather*} H_{\bb A^1}^{p,q}(X,Y,\bb Z)\cong\sheff\cc
O_{K^{Gr}}(M_{K^{Gr}}(X),M_{\cc K_0^{\oplus}}(Y)(q)[p-2q])\cong\\
\sheff\cc O_{K^\oplus}(M_{K^\oplus}(X),M_{\cc
K_0^{\oplus}}(Y)(q)[p-2q]).\end{gather*}
\end{lem}

\begin{proof}
All presheaves of the statement are semistable. The proof of
Theorem~\ref{nakonets} shows that the natural maps
   $$M_{K^{Gr}}(X)\to M_{K^\oplus}(X)\to M_{K}(X)$$
are isomorphisms in $\sheff(F)$. Now Corollary~\ref{vottak} and
Lemma~\ref{qwer} imply the claim.
\end{proof}

\begin{cor}\label{aqemm}
There is a natural isomorphism
   \begin{gather*} K_i(X)\cong\sheff\cc
O_{K^{Gr}}(M_{K^{Gr}}(X)[i],M_{K^{Gr}}(pt))\cong\\ \sheff\cc
O_{K^\oplus}(M_{K^\oplus}(X)[i],M_{K^\oplus}(pt))\cong\sheff\cc
O_{K}(M_{K}(X)[i],M_{K}(pt))
   \end{gather*}
for any $i\in\bb Z$ and $X\in Sm/F$.
\end{cor}

\begin{proof}
This follows from Corollary~\ref{trew} and Lemma~\ref{aqe}.
\end{proof}

\begin{cor}\label{qmq}
The maps of spectral categories
   $$\cc O_{K^{Gr}}\to\cc O_{K^{\oplus}}\to\cc O_K$$
induce triangulated equivalences
   $$\sheff\cc O_{K^{Gr}}\to\sheff\cc O_{K^{\oplus}}\to\sheff\cc O_K$$
of compactly generated triangulated categories.
\end{cor}

\begin{proof}
The objects $\{M_{K^{Gr}}(X)[i]\}_{i\in\bb Z,X\in Sm/F}$
(respectively $\{M_{K^\oplus}(X)[i]\}_{i\in\bb Z,X\in Sm/F}$ and
$\{M_{K}(X)[i]\}_{i\in\bb Z,X\in Sm/F}$) are compact generators of
the compactly generated triangulated category $\sheff\cc O_{K^{Gr}}$
(respectively $\sheff\cc O_{K^\oplus}$ and $\sheff\cc O_K$). Both
functors take the compact generators to compact generators and
induce isomorphisms of Hom-sets between them by Lemma~\ref{aqe}. Now
our assertion follows from standard facts about compactly generated
triangulated categories.
\end{proof}

We now have all the necessary information to prove the following
result saying that the Grayson (bivariant) motivic spectral sequence
is realized in a natural way in the triangulated category of
$K^{Gr}$-motives.

\begin{thm}\label{spektralka3}
For every $q\geq 0$ and every $X\in Sm/F$ the sequence of maps
   $$M_{K^{Gr}}(X)(q+1)\xrightarrow{f_{q+1}}M_{K^{Gr}}(X)(q)\lra{g_{q}}M_{\cc K_0^{\oplus}}(q)\lra{+}$$
is a triangle in $\sheff\cc O_{K^{Gr}}$ and the Grayson tower in
$\sheff\cc O_{K^{Gr}}$
   $$\cdots\xrightarrow{f_{q+2}}M_{K^{Gr}}(X)(q+1)\xrightarrow{f_{q+1}}M_{K^{Gr}}(X)(q)\lra{f_{q}}\cdots\lra{f_1}M_{K^{Gr}}(X)(0)$$
produces a strongly convergent spectral sequence
   $$E^2_{pq}=\sheff\cc O_{K^{Gr}}(M_{K^{Gr}}(U)[p+q],M_{\cc K_0^{\oplus}}(X)(q))\Rightarrow \sheff\cc O_{K^{Gr}}(M_{K^{Gr}}(U)[p+q],M_{K^{Gr}}(X)).$$
It agrees with the $\bb A^1$-local bivariant motivic spectral
sequence of Theorem~\ref{spektralka2}. Assume further that one of
the following conditions is satisfied:
\begin{enumerate}
\item $X=pt$;

\item the field $F$ is perfect.
\end{enumerate}
Then this spectral sequence agrees with the bivariant motivic
spectral sequence of Theorem~\ref{spektralka}.
\end{thm}

\begin{proof}
This follows from Theorem~\ref{spektralka2} and Lemma~\ref{aqe}.
\end{proof}

\section{Concluding remarks}

The interested reader may have observed that the authors have not
considered monoidal structures on the category of $\cc
O_{K^\oplus}$-modules. We believe that there should exist new
transfers $Cor_{virt}$ on $Sm/F$ which produce a spectral category
$\cc O_{K^{new}}$ such that:

\begin{itemize}
\item[$\diamond$] $\cc O_{K^{new}}$ is symmetric monoidal and
motivically excisive;

\item[$\diamond$] for any $X\in Sm/F$, $\cc O_{K^{new}}(-,X)$ is a
sheaf of symmetric spectra;

\item[$\diamond$] the motivic model category of (pre-)sheaves of symmetric spectra which are also
$\cc O_{K^{new}}$-modules is zig-zag Quillen equivalent to the
motivic model category for $\cc O_{K^\oplus}$-modules.
\end{itemize}

The use of motivically excisive spectral categories on $Sm/F$ and
their modules is a reminiscence of the theory of sheaves of $\cc
O_X$-modules over a ringed space $(X,\cc O_X)$. The structure sheaf
$\cc O_X$ is replaced with a motivically excisive spectral category
$\cc O$ on $Sm/F$ and sheaves of $\cc O_X$-modules are replaced with
the motivic model category of (pre-)sheaves of symmetric spectra
which are also $\cc O$-modules.

From this point of view the theory of spectral categories over
$Sm/F$ and their modules is a sort of ``motivic brave new algebra",
where the base symmetric monoidal model category is
$Pre^\Sigma(Sm/F)$.

\appendix
\section{Eilenberg--Mac Lane spectral categories}

In this section we construct the Eilenberg--Mac Lane spectral
categories $EM(\cc A)$ associated with ringoids $\cc A$ which are
equivalent to $H\cc A$. Although the authors have not found such
constructions in the literature, they do not have pretensions to
originality.

Let $(A,+)$ be an abelian monoid with neutral element 0 and let
$\Ar[n]$ be the category of arrows for the poset $\{0<1<\cdots<n\}$.
It can be regarded as the partially ordered set of pairs $(i,j)$,
$0\leq i\leq j\leq n$, where $(i,j)\leq(i',j')$ \ifff $i\leq i'$ and
$j\leq j'$.

We consider the set of functions
   \begin{align*}
    a:\Ob\Ar[n]&\lra{}A\\
    (i,j)&\longmapsto a(i,j)=a_{i,j}
   \end{align*}
having the property that for every $j$, $a_{j,j}=0$ and
$a_{i,k}=a_{i,j}+a_{j,k}$ whenever $i\leq j\leq k$. Let us denote it
by $\sigma_nA$.

Given a map $u:[m]\lra{}[n]$ in $\Delta$, the function
$u^*:\sigma_nA\lra{}\sigma_mA$ sends the element $(i,j)\mapsto
a(i,j)$ to the element $(r,s)\longmapsto a(u(r),u(s))$.

The elements of $\sigma_nA$ may also be regarded as diagrams of the
form
   \begin{equation}\label{rrr}
    \begin{matrix}
     {}&{}&{}&{}&{}&{}&{}&{}&a_{n-1,n}\\
     {}&{}&{}&{}&{}&{}&{}&{}&{}\\
     {}&{}&{}&{}&{}&{}&{}&{}&\vdots\\
     {}&{}&{}&{}&a_{23}&{}&\cdots&{}&a_{2n}\\
     {}&{}&&&{}&&&&{}\\
     {}&{}&a_{12}&{}&a_{13}&{}&\cdots&{}&a_{1n}\\
     {}&{}&{}&&{}&&&&{}\\
     a_{01}&{}&a_{02}&{}&a_{03}&{}&\cdots&{}&a_{0n}
    \end{matrix}
   \end{equation}
Then the degeneracy maps are defined as the functions
$s_i:\sigma_nA\to\sigma_{n+1}A$ by duplicating $a_{0,i}$, and
reindexing with the normalization $a_{i,i+1}=0$.

Also, the face map $d_0:\sigma_nA\to\sigma_{n-1}A$ is the function
which is defined by deleting the bottom row of~\eqref{rrr}. For
$0<i\le n$ we define the face maps as the functions
$d_i:\sigma_nA\to\sigma_{n-1}A$ by omitting the row $a_{i,*}$ and
the column containing $a_{0i}$ in~\eqref{rrr}, and reindexing the
$a_{j,k}$ as needed.

Given two functions $a,b:\Ob\Ar[n]\to A$, we define a binary
operation on $\sigma_nA$ as $(a+b)_{i,j}:=a_{i,j}+b_{i,j}$. Then the
set $\sigma_nA$ is an abelian monoid as well. So we arrive at a
simplicial abelian monoid $n\mapsto\sigma_nA$ denoted by $\sigma.A$.
We can iterate the $\sigma$-construction to get a bisimplicial
abelian monoid $\sigma^2A=\sigma.\sigma.A$ or, more generally, a
multisimplicial abelian monoid $\sigma^nA$.

Similar to Waldhausen's $S$-construction~\cite{Wal} there is a
natural inclusion $A\wedge S^1\to\sigma.A$, and by adjointness
therefore an inclusion of $|A|$ into the loop space of $|\sigma.A|$.
There results a spectrum
   $$EM(A):=(A,\sigma.A,\sigma.\sigma.A,\ldots)$$
whose structural maps are defined just as the map
$|A|\to\Omega|\sigma.A|$ above. One can actually define $EM(A)$ as a
symmetric spectrum in a similar way that in~\cite{GH}.

Let $B,C$ be two other abelian monoids. A map $f:A\times B\to C$ is
called a {\it bilinear pairing\/} if $f(a+a',b)=f(a,b)+f(a',b)$ and
$f(a,b+b')=f(a,b)+f(a,b')$ for all $a,a'\in A$ and $b,b'\in B$. In
particular, $f(a,0)=f(0,b)=0$. Every bilinear map induces a map of
symmetric spectra
   $$f:EM(A)\wedge EM(B)\to EM(C).$$
These maps are associative for strictly associative bilinear
pairings.

The universal example of an abelian monoid which acts on any other
such monoid is the monoid of non-negative integers $\bb Z_{\geq 0}$.
Given an abelian monoid $A$, there is a bilinear pairing $\bb
Z_{\geq 0}\times A\to A$ sending $(n,a)$ to $na$.

Recall that a {\it semiring\/} is a set $A$ equipped with two binary
operations + and $*$, called addition and multiplication, such that
$(A,+)$ is an abelian monoid with identity element 0, $(A,*)$ is a
monoid with identity element 1, multiplication distributes over
addition, and 0 annihilates $A$ with respect to multiplication. It
follows from above that $EM(A)$ is a symmetric ring spectrum for
every semiring $A$ such that the structure map from the sphere
spectrum to $EM(A)$ is given by the map $S^0\to EM(A)_0=A$ sending
the basepoint to 0 and the non-basepoint to 1.

The universal example of a semiring which is mapped to any other
semiring is $\bb Z_{\geq 0}$. One easily sees that for every abelian
monoid $A$ the symmetric spectrum $EM(A)$ is an $EM(\bb Z_{\geq
0})$-module. If $A$ is a semiring then there is a natural map of
ring spectra
   $$EM(\bb Z_{\geq 0})\to EM(A).$$

More generally, any {\it semiringoid\/} $\cc A$, that is a category
whose Hom-sets are abelian monoids with bilinear composition and
whose End-sets are semirings gives rise to a spectral category
$EM(\cc A)$ with $EM(\cc A)(x,y)=EM(\Hom_{\cc A}(x,y))$ for all
$x,y\in\Ob\cc A$.

Let $i\cc C$ be the Waldhasen category of isomorphisms for a
Waldhausen category $\cc C$ in which all cofibrations split. For
instance, $\cc C$ is an additive category or the category $\Gamma$
of finite pointed sets $n^+=\{0,1,\ldots,n\}$ with 0 as basepoint
and pointed set maps. Denote by $\pi_0\cc C$ the abelian monoid of
isomorphism classes for objects in $\cc C$ (e.g., $\pi_0\Gamma=\bb
Z_{\geq 0}$). Given two classes $[A],[B]\in\pi_0\cc C$, the binary
operation is defined as usual $[A]+[B]:=[A\coprod B]$. There is a
natural map of spectra
   $$\tau:K(\cc C)\to EM(\pi_0\cc C)$$
sending each diagram $(F:\Ar[n_1]\times\cdots\times\Ar[n_d]\to\cc
C)\in S^d\cc C$ to the composition
$(\Ob\Ar[n_1]\times\cdots\times\Ob\Ar[n_d]\bl F\to\Ob\cc
C\to\pi_0\cc C)\in\sigma^d(\pi_0\cc C)$.

Let $i\cc A,i\cc B$ be other two Waldhasen categories of
isomorphisms for Waldhausen categories $\cc A,\cc B$ in which all
cofibrations split. Suppose $f:\cc A\times\cc B\to\cc C$ is a
biexact functor between Waldhausen categories. Then the following
diagram of maps of symmetric spectra commutes
   $$\xymatrix{K(\cc A)\wedge K(\cc B)\ar[r]^(.65)f\ar[d]_{\tau\wedge\tau}&K(\cc C)\ar[d]^\tau\\
               EM(\pi_0\cc A)\wedge EM(\pi_0\cc B)\ar[r]^(.65)f&EM(\pi_0\cc C).}$$

Let $\cc O$ be a spectral category generated by Waldhausen
categories of isomorphisms $i\cc C(x,y)$, where $x,y\in\Ob\cc O$, in
which all cofibrations split (i.e., $\cc O(x,y)=K(\cc C(x,y))$), and
such that the composition law
   $$\cc O(y,z)\wedge\cc O(x,y)\to\cc O(x,z)$$
comes from biexact functors $\cc C(y,z)\times\cc C(x,y)\to\cc
C(x,z)$. For example, $\cc O=\cc O_{K^\oplus}$. Then the collections
$\pi_0(\cc C(x,y))$, $x,y\in\Ob\cc O$, are abelian monoids and form
a semiringoid. It gives rise to a spectral category, denoted by
$EM(\pi_0\cc O)$. There is a map of spectral categories
   $$\upsilon:\cc O\to EM(\pi_0\cc O).$$

The additivity theorem for the $\sigma$-construction says that the
natural map
   $$\sigma.(\sigma_2A)\to\sigma.A\times\sigma.A$$
sending $a:\Ob\Ar[2]\to A$ to the couple $(a_{0,1},a_{1,2})$ is an
isomorphism of simplicial sets for every abelian monoid $A$.
Repeating now Waldhausen's results~\cite[section~1.5]{Wal} for the
$\sigma.$-construction, we get that the natural map
$|\sigma.A|\to\Omega|\sigma.\sigma.A|$ which is adjoint to the map
$|\sigma.A|\wedge S^1\to |\sigma.\sigma.A|$ is a homotopy
equivalence and more generally also the map
$|\sigma^nA|\to\Omega|\sigma^{n+1}A|$ for every $n>0$. This proves
that the spectrum $EM(A)$ is an $\Omega$-spectrum beyond the first
term.

There is an isomorphism of simplicial sets $\sigma.A\cong BA$, where
$BA$ stands for the classifying space of $A$. It takes every
diagram~\eqref{rrr} to $(a_{0,1},a_{1,2},\ldots,a_{n-1,n})\in BA$.
We conclude that if $A$ is an abelian group, then each space
$|\sigma^nA|$, $n\geq 0$, has the homotopy type of the
Eilenberg--Mac~Lane space $K(A,n)$. Moreover, $EM(A)$ is a genuine
$\Omega$-spectrum.

Consider a spectral category $\cc O$ generated by Waldhausen
categories (see above). Then the collections $K_0(\cc
C(x,y))=\pi_0(\cc O(x,y))$, $x,y\in\Ob\cc O$, are abelian groups and
form a ringoid. It gives rise to a spectral category, denoted by
$K_0\cc O$. There are maps of spectral categories
   $$\cc O\bl\upsilon\to EM(\pi_0\cc O)\xrightarrow\kappa K_0\cc O,$$
where $\kappa$ is induced by the universal group completion maps
$\pi_0(|i\cc O(x,y)|)\to K_0(\cc C(x,y))$.

Consider the symmetric spectrum $HA$ associated with an abelian
group $A$. Recall that $HA_p=A\otimes\wt{\bb Z}[S^p]$ for any $p\geq
0$. The identity map of $A$ induces a map of symmetric spectra
   $$\l:\Sigma^\infty A=(A,A\wedge S^1,A\wedge S^2,\ldots)\to EM(A).$$
The maps $l_p:A\wedge S^p\to|\sigma^pA|$ induce in a unique way maps
$\ell_p:A\otimes\wt{\bb Z}[S^p]\to|\sigma^pA|$. These yield a map of
symmetric spectra
   $$\ell_A:HA\to EM(A),$$
functorial in $A$. Recall that $HA_p$ has homotopy type of the
Eilenberg--Mac~Lane space $K(A,p)$. We deduce that each $\ell_p$ is
a weak equivalence, and hence $\ell$ is a levelwise weak
equivalence.

For two abelian groups $A$ and $B$, there is a natural morphism of
symmetric spectra (see, e.g.,~\cite[Example~I.3.11]{Sch})
   $$HA\wedge HB\to H(A\otimes B).$$
It is easily verified that the diagram
   $$\xymatrix{HA\wedge HB\ar[r]\ar[d]_{\ell_A\wedge\ell_B}&H(A\otimes B)\ar[d]^{\ell_{A\otimes B}}\\
               EM(A)\wedge EM(B)\ar[r]&EM(A\otimes B).}$$
is commutative.

Now let $\cc A$ be a ringoid and let $H\cc A$ be the spectral
category associated with it. Recall that $H\cc A(x,y)_p={\cc
A}(x,y)\otimes\wt{\bb Z}[S^p]$ for any objects $x,y\in\Ob\cc A$ and
$p\geq 0$. It follows from above that there is a map of spectral
categories
   $$\ell:H\cc A\to EM(\cc A)$$
such that $H\cc A(x,y)\to EM(\cc A)(x,y)$ is a levelwise equivalence
of symmetric spectra for any $x,y\in\Ob\cc A$.

\end{document}